\documentclass[reqno,a4paper,11pt]{article}
 \usepackage{amssymb,amsmath}
 \usepackage{amsfonts}
 \usepackage{mathrsfs,enumerate}
 \usepackage{color}
 \usepackage[
 color,
 final] 
 {showkeys}

 \definecolor{refkeybis}{gray}{.65}
 \definecolor{labelkeybis}{gray}{.65}
 {\makeatletter
 \def\SK@refcolor{\color{refkeybis}}%
 \def\SK@labelcolor{\color{labelkeybis}}}

 \numberwithin{equation}{section} 
 \newtheorem{theorem}{Theorem}[section] 
 \newtheorem{lemma}[theorem]{Lemma}
 \newtheorem{corollary}[theorem]{Corollary}
 \newtheorem{proposition}[theorem]{Proposition}
 
 \newtheorem{remark}[theorem]{Remark}

 \newtheorem{definition}[theorem]{Definition}

 \newcommand{\R}{\mathbb{R}}
 \newcommand{\T}{\mathbb{T}}

 \newcommand{\N}{\mathbb{N}}
 \newcommand{\Z}{\mathbb{Z}}

 \renewcommand{\AA}{\mathscr{A}}

 \newcommand{\MM}{\mathscr{M}}

 \newcommand{\ii}{{\mbox{\boldmath$i$}}}
 \newcommand{\mm}{{\mbox{\boldmath$m$}}}

 \newcommand{\sii}{{\mbox{\scriptsize\boldmath$i$}}}
 
 \newcommand{\uu}{{\mbox{\boldmath$u$}}}
 \newcommand{\vv}{{\mbox{\boldmath$v$}}}
 \newcommand{\ww}{{\mbox{\boldmath$w$}}}
 \newcommand{\svv}{{\mbox{\scriptsize\boldmath$v$}}}
 
 \newcommand{\sww}{{\mbox{\scriptsize\boldmath$w$}}}

 \newcommand{\tauV}{{\kern-3pt\tau}}

 \newcommand{\oVVVk}{\overline{\mbox{\boldmath$V$}}\kern-3pt}
 \newcommand{\tVVVk}{\tilde{\mbox{\boldmath$V$}}\kern-3pt}

 \newcommand{\llambda}{{\mbox{\boldmath$\lambda$}}}

 \newcommand{\mmu}{{\mbox{\boldmath$\mu$}}}
 \newcommand{\nnu}{{\mbox{\boldmath$\nu$}}}
 \newcommand{\snnu}{{\mbox{\scriptsize\boldmath$\nu$}}}

 \newcommand{\seeta}{{\mbox{\scriptsize\boldmath$\eta$}}}

 \newcommand{\eeta}{{\mbox{\boldmath$\eta$}}}

 \renewcommand{\restriction}[1]{\lower3pt\hbox{$|_{#1}$}}
 
 \newcommand{\Leb}[1]{{\mathscr L}^{#1}} 
 \newcommand{\Haus}[1]{{\mathscr H}^{#1}} 

 \renewcommand{\div}{{\rm div}}
 \newcommand{\eps}{{\varepsilon}}
 \renewcommand{\to}{\rightarrow}
 \newenvironment{Proof}{\removelastskip\par\medskip 
 \noindent{\em Proof.}
 \rm}{\penalty-20\null\hfill$\square$\par\medbreak} 

 \renewcommand{\a}{\alpha}
 
 \renewcommand{\d}{\delta}
 
 \newcommand{\e}{\varepsilon}
 \newcommand{\g}{\gamma}
 \newcommand{\G}{\Gamma}
 \renewcommand{\l}{\lambda}
 
 \newcommand{\n}{\nabla}
 
 \newcommand{\s}{\sigma}

 \renewcommand{\O}{\Omega}
 \renewcommand{\o}{\omega}
 
 \renewcommand{\i}{\infty}
 \newcommand{\p}{\partial}
 \newcommand{\ov}{\overline}

 \newcommand{\dist}{\operatorname{dist}}
 \newcommand{\supp}{\operatorname{supp}}
 
 \newcommand{\BorelSets}[1]{\mathcal B(#1)}
 \newcommand{\Probabilities}[1]{\mathscr P(#1)} 
 \newcommand{\Measures}[1]{\mathscr M(#1)} 
 \newcommand{\PlusMeasures}[1]{\mathscr M_+(#1)} 
 
 \newcommand{\SDiff}[1]{{\rm SDiff}(#1)}
 \newcommand{\Oa}{\tilde\Omega}
 \makeindex
 
\begin{document}

 \title{Geodesics in the space of
 measure-preserving maps and plans
 \thanks{This work was partially supported by the MIUR COFIN04 grant}
 }
 \author{Luigi Ambrosio\
   \thanks{\textsf{l.ambrosio@sns.it}}\\
   Scuola Normale Superiore\\
   Pisa
   \and
   Alessio Figalli\
   \thanks{\textsf{a.figalli@sns.it}}\\
   Scuola Normale Superiore\\
   Pisa}
 \date{January 22, 2007}

 \maketitle

 \begin{abstract}
 We study Brenier's variational models for incompressible Euler equations.
 These models give rise to a relaxation of the Arnold distance in the space
 of measure-preserving maps and, more generally, measure-preserving plans.
 We analyze the properties of the relaxed distance, we
 show a close link between the Lagrangian and the Eulerian model, and we derive
 necessary and sufficient optimality conditions for minimizers. These
 conditions take into account a modified Lagrangian induced by the pressure field.
 Moreover, adapting some ideas of Shnirelman, we show that, even for
 non-deterministic final conditions, generalized
 flows can be approximated in energy by flows associated to measure-preserving
 maps.
 \end{abstract}


 \section{Introduction}

 The velocity of an incompressible fluid moving inside a region $D$
 is mathematically described by a time-dependent and
 divergence-free vector field $\uu(t,x)$ which is parallel to the
 boundary $\p D$. The Euler equations for incompressible fluids
 describes the evolution of such velocity field $\uu$ in terms of
 the pressure field $p$:
 \begin{equation}\label{Eu}
 \left\{
 \begin{array}{ll}
 \p_t \uu + (\uu \cdot \n)\uu= - \n p & \text{in }[0,T] \times D,\\
 \div\,\uu=0 & \text{in }[0,T] \times D,\\
 \uu \cdot n=0 & \text{on }[0,T] \times \p D.
 \end{array}
 \right.
 \end{equation}
 Let us assume that $\uu$ is smooth, so that it produces a unique
 flow $g$, given by
 $$
 \left\{
 \begin{array}{l}
 \dot g(t,a)=\uu(t,g(t,a)),\\
 g(0,a)=a.
 \end{array}
 \right.
 $$
 By the incompressibility condition, we get that at each time $t$
 the map $g(t,\cdot):D \to D$ is a measure-preserving
 diffeomorphism of $D$, that is
 $$
 g(t,\cdot)_\#\mu_D=\mu_D,
 $$
 (here and in the sequel $f_\#\mu$ is the push-forward of a measure
 $\mu$ through a map $f$, and $\mu_D$ is the volume measure of the
 manifold $D$). Writing Euler equations in terms of $g$, we get
 \begin{equation}\label{Eug}
 \left\{
 \begin{array}{ll}
 \ddot g(t,a) =-\n p\left(t,g(t,a)\right) & \text{$(t,a)\in [0,T] \times D$,}\\
 g(0,a)=a &\text{$a\in D$,}\\
 g(t,\cdot)\in \SDiff{D} &\text{$t\in [0,T]$.}
 \end{array}
 \right.
 \end{equation}
 Viewing the space $\SDiff{D}$ of measure-preserving
 diffeomorphisms of $D$ as an infinite-dimensional manifold with
 the metric inherited from the embedding in $L^2$, and with tangent
 space made by the divergence-free vector fields, Arnold
 interpreted the equation above, and therefore \eqref{Eu}, as a
 \emph{geodesic} equation on $\SDiff{D}$ \cite{arnold}. According
 to this intepretation, one can look for solutions of \eqref{Eug}
 by minimizing
 \begin{equation}\label{Eugg}
 T\int_0^T\int_D \frac{1}{2} |\dot g(t,x)|^2\,d\mu_D(x)\,dt
 \end{equation}
 among all paths $g(t,\cdot):[0,T]\to\SDiff{D}$ with $g(0,\cdot)=f$
 and $g(T,\cdot)=h$ prescribed (typically, by right invariance, $f$
 is taken as the identity map $\ii$), and the pressure field arises
 as a Lagrange multiplier from the incompressibility constraint
 (the factor $T$ in front of the integral is just to make the
 functional scale invariant in time). We shall denote by
 $\delta(f,h)$ the Arnold distance in $\SDiff{D}$, whose
 square is defined by the
 above-mentioned variational problem in the time interval $[0,1]$.

 Although in the traditional approach to \eqref{Eu} the initial
 velocity is prescribed, while in the minimization of \eqref{Eugg}
 is not, this variational problem has an independent interest and
 leads to deep mathematical questions, namely existence of relaxed
 solutions, gap phenomena and necessary and sufficient optimality
 conditions, that are investigated in this paper. We also remark
 that \emph{no} existence result of distributional solutions of
 \eqref{Eu} is known when $d>2$ (the case $d=2$ is different,
 thanks to the vorticity formulation of \eqref{Eu}), see
 \cite{lions}, \cite{brenierextra} for a discussion on this topic
 and other concepts of weak solutions to \eqref{Eu}.

 On the positive side, Ebin and Marsden proved in
 \cite{ebinmarsden} that, when $D$ is a smooth compact manifold
 with no boundary, the minimization of \eqref{Eugg} leads to a
 unique solution, corresponding also to a solution to Euler
 equations, if $f$ and $h$ are sufficienly close in a suitable
 Sobolev norm.

 On the negative side, Shnirelman proved in \cite{shnir1},
 \cite{shnir2} that when $d\geq 3$ the infimum is not attained in
 general, and that when $d=2$ there exists $h\in\SDiff{D}$ which
 cannot be connected to $\ii$ by a path with finite action. These
 ``negative'' results motivate the study of relaxed versions of
 Arnold's problem.

 The first relaxed version of Arnold's minimization problem was
 introduced by Brenier in \cite{brenier1}: he considered
 probability measures $\eeta$ in $\O(D)$, the space of continuous
 paths $\omega:[0,T]\to D$, and minimized the energy
 $$
 \AA_T(\eeta):=T\int_{\O(D)}\int_0^T\frac{1}{2}|\dot\omega(\tau)|^2\,d\tau\,d\eeta(\omega),
 $$
 with the constraints
 \begin{equation}\label{junior}
 (e_0,e_T)_\#\eeta=(\ii,h)_\#\mu_D,\qquad
 (e_t)_\#\eeta=\mu_D\quad\forall t\in [0,T]
 \end{equation}
 (here and in the sequel $e_t(\omega):=\omega(t)$ are the
 evaluation maps at time $t$). According to Brenier, we shall call
 these $\eeta$ \emph{generalized incompressible flows} in $[0,T]$
 between $\ii$ and $h$. Obviously any sufficiently regular path
 $g(t,\cdot):[0,1]\to S(D)$ induces a generalized incompressible
 flow $\eeta=(\Phi_g)_\#\mu_D$, where $\Phi_g:D\to\O(D)$ is given
 by $\Phi_g(x)=g(\cdot,x)$, but the converse is far from being
 true: the main difference between classical and generalized flows
 consists in the fact that fluid paths starting from different
 points are allowed to cross at a later time, and fluid paths
 starting from the same point are allowed to split at a later time.
 This approach is by now quite common, see for instance \cite{ambrosio}
 (DiPerna-Lions theory), \cite{morel} (branched optimal transportation),
 \cite{lottvillani}, \cite{villaniSF}.

 Brenier's formulation makes sense not only if $h\in\SDiff{D}$, but
 also when $h\in S(D)$, where $S(D)$ is the space of
 measure-preserving maps $h:D\to D$, not necessarily invertible or
 smooth. In the case $D=[0,1]^d$, existence of admissible paths
 with finite action connecting $\ii$ to any $h\in S(D)$ was proved
 in \cite{brenier1}, together with the existence of paths with minimal action.
 Furthermore, a consistency result was proved: smooth solutions to
 \eqref{Eu} are optimal even in the larger class of the generalized
 incompressible flows, provided the pressure field $p$ satisfies
 \begin{equation}\label{goodp}
 T^2\sup_{t\in [0,T]}\sup_{x\in D}|\nabla^2_xp(t,x)| \leq \pi^2,
 \end{equation}
 and are the unique ones if the inequality is strict.
 When $\eeta=(\Phi_g)_\#\mu_D$ we can recover $g(t,\cdot)$ from
 $\eeta$ using the identity
 $$
 (e_0,e_t)_\#\eeta=(\ii,g(t,\cdot))_\#\mu_D, \qquad t\in[0,T].
 $$
 Brenier found in \cite{brenier1} examples of action-minimizing
 paths $\eeta$ (for instance in the unit ball of $\R^2$, between
 $\ii$ and $-\ii$) where no such representation is possible. The
 same examples show that the upper bound \eqref{goodp} is sharp.
 Notice however that $(e_0,e_t)_\#\eeta$ is a measure-preserving
 plan, i.e. a probability measure in $D\times D$ having both
 marginals equal to $\mu_D$. Denoting by $\Gamma(D)$ the space of
 measure-preserving plans, it is therefore natural to consider
 $t\mapsto (e_0,e_t)_\#\eeta$ as a ``minimizing geodesic'' between
 $\ii$ and $h$ in the larger space of measure-preserving plans.
 Then, to be consistent, one has to extend Brenier's minimization
 problem considering paths connecting $\gamma,\,\eta\in\Gamma(D)$.
 We define this extension, that reveals to be useful also to
 connect this model to the Eulerian-Lagrangian one in
 \cite{brenier4}, and to obtain necessary and sufficient optimality
 conditions even when only ``deterministic'' data $\ii$ and $h$ are
 considered (because, as we said, the path might be
 non-deterministic in between). In this presentation of our
 results, however, to simplify the matter as much as possible, we
 shall consider the case of paths $\eeta$ between $\ii$ and $h\in
 S(D)$ only.

 In Section~\ref{sgap} we study the relation between the relaxation
 $\delta_*$ of the Arnold distance, defined by
 $$
 \delta_*(h):=\inf\left\{\liminf_{n\to\infty}\delta(\ii,h_n):\
 h_n\in\SDiff{D},\,\,\int_D|h_n-h|^2\,d\mu_D\to 0\right\},
 $$
 and the distance $\overline\delta(\ii,h)$ arising from the
 minimization of the Lagrangian model. It is not hard to show that
 $\overline\delta(\ii,h)\leq\delta_*(h)$, and a natural question is
 whether equality holds, or a gap phenomenon occurs. In the case
 $D=[0,1]^d$ with $d>2$, an important step forward was obtained by
 Shnirelman in \cite{shnir2}, who proved that equality holds when
 $h\in\SDiff{D}$; Shnirelman's construction provides an
 approximation (with convergence of the action) of generalized
 flows connecting $\ii$ to $h$ by smooth flows still connecting
 $\ii$ to $h$. The main result of this section is the proof that no
 gap phenomenon occurs, still in the case $D=[0,1]^d$ with $d>2$,
 even when non-deterministic final data (i.e. measure-preserving
 plans) are considered. The proof of this fact is based on an
 auxiliary approximation result, Theorem~\ref{nogapGamma2},
 valid in any number of dimensions, which we believe of independent
 interest: it allows to approximate, with convergence of the
 action, any generalized flow $\eeta$ in $[0,1]^d$ by $W^{1,2}$
 flows (in time) induced by measure-preserving maps $g(t,\cdot)$.
 This fact shows that the ``negative'' result of Shnirelman on the existence in
 dimension $2$ of non-attainable diffeomorphisms is due to the regularity assumption on the path,
 and it is false if one allows for paths in the larger space $S(D)$.
 The proof of Theorem~\ref{nogapGamma2} uses some key ideas from \cite{shnir2} (in particular
 the combination of law of large numbers and smoothing of discrete
 families of trajectories), and some ideas coming from the theory
 of optimal transportation.

 Minimizing generalized paths $\eeta$ are not unique in general, as
 shown in \cite{brenier1}; however, Brenier proved in
 \cite{brenier2} that the gradient of the pressure field $p$, identified by the
 distributional relation
 \begin{equation}\label{mezzetti}
 \nabla
 p(t,x)=-\partial_t\overline{\vv}_t(x)-\div\left(\overline{\vv\otimes\vv}_t(x)\right),
 \end{equation}
 is indeed unique. Here $\overline{\vv}_t(x)$ is the ``effective
 velocity'', defined by $(e_t)_\#(\dot\omega(t)\eeta)=\overline{\vv}_t\mu_D$, and
 $\overline{\vv\otimes\vv}_t$ is the quadratic effective velocity,
 defined by $(e_t)_\#(\dot\omega(t)\otimes\dot\omega(t)\eeta)=\overline{\vv\otimes\vv}_t\mu_D$.
 The proof of this fact is based on the so-called dual least action
 principle: if $\eeta$ is optimal, we have
 \begin{equation}\label{dualp}
 \AA_T(\nnu)\geq\AA_T(\eeta)+\langle p,\rho^{\snnu}-1\rangle
 \end{equation}
 for any measure $\nnu$ in $\O(D)$ such that
 $(e_0,e_T)_\#\nnu=(\ii,h)_\#\mu_D$ and
 $\Vert\rho^{\snnu}-1\Vert_{C^1}\leq 1/2$. Here $\rho^{\snnu}$ is
 the (absolutely continuous) density produced by the flow $\nnu$,
 defined by $\rho^{\snnu}(t,\cdot)\mu_D=(e_t)_\#\nnu$. In this way,
 the incompressibility constraint can be slightly relaxed and one
 can work with the augmented functional (still minimized by
 $\eeta$)
 $$
 \nnu\mapsto\AA_T(\nnu)-\langle p,\rho^{\snnu}-1\rangle,
 $$
 whose first variation leads to \eqref{mezzetti}.

 In Theorem~\ref{thmslighcomprflows}, still using the key
 Proposition~2.1 from \cite{brenier2}, we provide a simpler proof
 and a new interpretation of the dual least action principle.

 A few years later, Brenier introduced in \cite{brenier4} a new
 relaxed version of Arnold's problem of a mixed Eulerian-Lagrangian
 nature: the idea is to add to the Eulerian variable $x$ a
 Lagrangian one $a$ representing, at least when $f=\ii$, the
 initial position of the particle; then, one minimizes a functional
 of the Eulerian variables (density and velocity), depending also on $a$.
 Brenier's motivation for looking at the new model was that this
 formalism allows to show much stronger regularity results for the
 pressure field, namely $\partial_{x_i} p$ are locally finite measures
 in $(0,T)\times D$.
 In Section~\ref{secConstr} we describe in detail this new model and, in
 Section~\ref{secEquivalence}, we show that the two models are
 basically equivalent. This result will be used by us to transfer
 the regularity informations on the pressure field up to the
 Lagrangian model, thus obtaining the validity of \eqref{dualp} for
 a much larger class of generalized flows $\nnu$, that we call
 flows with bounded compression. The proof of the equivalence
 follows by a general principle (Theorem~\ref{tsuperp}, borrowed
 from \cite{amgisa}) that allows to move from an Eulerian to a
 Lagrangian description, lifting solutions to the continuity
 equation to measures in the space of continuous maps.

 In the final section of our paper we look for necessary and
 sufficient optimality conditions for the geodesic problem. These
 conditions require that the pressure field $p$ is a function and
 not only a distribution: this technical result is achieved in \cite{amfiga},
 where, by carefully analyzing and improving Brenier's difference-quotient argument,
 we show that $\partial_{x_i} p\in L^2_{\rm loc}\bigl((0,T);{\cal M}_{\rm loc}(D)\bigr)$
 (this implies, by Sobolev embedding, $p\in L^2_{\rm loc}\bigl((0,T);L^{d/(d-1)}_{\rm loc}(D)\bigr)$).

 In this final section, although we do not see a serious obstruction to the extension of
 our results to a more general framework, we consider the case of the flat torus $\T^d$ only, and we shall denote by
 $\mu_\T$ the canonical measure on the flat torus.
 We observe that in this case $p\in L^2_{\rm loc}\bigl((0,T);L^{d/(d-1)}(\T^d)\bigr)$ and so, taking
 into account that the pressure field in (\ref{dualp}) is uniquely determined up to additive time-dependent constants,
 we may assume that $\int_{\T^d} p(t,\cdot)d\mu_\T=0$ for almost all $t \in (0,T)$.\\
 The first elementary remark is that any integrable function $q$ in
 $(0,T)\times\T^d$ with $\int_{\T^d}q(t,\cdot)\,d\mu_\T=0$ for
 almost all $t\in (0,T)$ provides us with a null-lagrangian for the
 geodesic problem, as the incompressibility constraint gives
 $$
 \int_{\O(\T^d)}\int_0^T q(t,\omega(t))\,dt\,d\nnu(\omega)=
 \int_0^T \int_{\T^d} q(t,x)\,d\mu_\T(x)\,dt=0
 $$
 for any generalized incompressible flow $\nnu$. Taking also the
 constraint $(e_0,e_T)_\#\nnu=(\ii,h)_\#\mu$ into account, we get
 $$
 \AA_T(\nnu)=T\int_{\O(\T^d)}\left(
 \int_0^T\frac{1}{2}|\dot\omega(t)|^2-q(t,\omega)\,dt\right)
 \,d\nnu(\omega)\geq \int_{\T^d}c^T_q(x,h(x))\,d\mu_\T(x),
 $$
 where $c^T_q(x,y)$ is the minimal cost associated with the
 Lagrangian $T\int_0^T\frac{1}{2}|\dot\omega(t)|^2-q(t,\omega)\,dt$. Since
 this lower bound depends only on $h$, we obtain that any $\eeta$
 satisfying \eqref{junior} and
 concentrated on $c_q$-minimal paths, for some $q\in L^1$, is
 optimal, and $\overline\delta^2(\ii,h)=\int c^T_q(\ii,h)\,d\mu_\T$.
 This is basically the argument used by Brenier in \cite{brenier1}
 to show the minimality of smooth solutions to \eqref{Eu}, under
 assumption \eqref{goodp}: indeed, this condition guarantees that
 solutions of $\ddot\omega(t)=-\nabla p(t,\omega)$ (i.e. stationary
 paths for the Lagrangian, with $q=p$) are also minimal.

 We are able to show that basically this condition is
 \emph{necessary and sufficient} for optimality if the pressure
 field is globally integrable (see Theorem~\ref{optima3}). However,
 since no global in time regularity result for the pressure field
 is presently known, we have also been looking for necessary and
 sufficient optimality conditions that don't require the global
 integrability of the pressure field. Using the regularity $p\in
 L^1_{\rm loc}\left((0,T);L^r(D)\right)$ for some $r>1$, guaranteed
 in the case $D=\T^d$ with $r=d/(d-1)$ by \cite{amfiga}, we show
 in Theorem~\ref{optima1} that any optimal $\eeta$ is concentrated
 on \emph{locally minimizing} paths for the Lagrangian
 \begin{equation}
 \label{defLp}
 {\cal L}_p(\omega):=\int
 \frac{1}{2}|\dot\omega(t)|^2-p(t,\omega)\,dt
 \end{equation}
 Since we are going to integrate $p$ along curves, this statement
 is not invariant under modifications of $p$ in negligible sets,
 and the choice of a \emph{specific} representative
 $\bar p(t,x):=\liminf_{\e\downarrow 0}p(t,\cdot)\ast\phi_\e(x)$
 in the Lebesgue
 equivalence class is needed. Moreover, the necessity of
 pointwise uniform estimates on $p_\e$ requires the integrability
 of $Mp(t,x)$, the maximal function of $p(t,\cdot)$ at $x$ (see \eqref{maximal}).

 In addition, we identify a second necessary (and more hidden)
 optimality condition. In order to state it, let us consider an
 interval $[s,t]\subset (0,T)$ and the cost function
 \begin{equation}\label{olla*}
 c^{s,t}_p(x,y):=\inf\left\{
 \int_s^t\frac{1}{2}|\dot\omega(\tau)|^2-p(\tau,\omega)\,d\tau:\
 \omega(s)=x,\,\omega(t)=y,\,Mp(\tau,\omega)\in L^1(s,t)\right\}.
 \end{equation}
 (the assumption $Mp(\tau,\omega)\in L^1(s,t)$ is forced by technical reasons).
 Recall that, according to the theory
 of optimal transportation, a probability measure $\lambda$ in
 $\T^d\times\T^d$ is said to be $c$-optimal if
 $$
 \int_{\T^d\times\T^d}c(x,y)\,d\lambda'\geq
 \int_{\T^d\times\T^d}c(x,y)\,d\lambda
 $$
 for any probability measure $\lambda'$ having the same marginals $\mu_1$, $\mu_2$
 of $\lambda$. We shall also denote $W_c(\mu_1,\mu_2)$ the minimal value, i.e.
 $\int_{\T^d\times\T^d}c\,d\lambda$, with $\lambda$ $c$-optimal.
 Now, let $\eeta$ be an optimal generalized
 incompressible flow between $\ii$ and $h$; according to the
 disintegration theorem, we can represent
 $\eeta=\int\eeta_a\,d\mu_D(a)$, with $\eeta_a$ concentrated on
 curves starting at $a$ (and ending, since our final conditions is
 deterministic, at $h(a)$), and consider the plans
 $\l^{s,t}_a=(e_s,e_t)_\#\eeta_a$. We show that
 \begin{equation}\label{edura}
 \text{for all $[s,t]\subset (0,T)$, \quad $\l^{s,t}_a$ is
 $c_p^{s,t}$-optimal for $\mu_\T$-a.e. $a\in\T^d$.}
 \end{equation}
 Roughly speaking, this condition tells us that one has not only to
 move mass from $x$ to $y$ achieving $c_p^{s,t}$, but also to
 optimize the distribution of mass between time $s$ and time $t$.
 In the ``deterministic'' case when either $(e_0,e_s)_\#\eeta$ or
 $(e_0,e_t)_\#\eeta$ are induced by a transport map $g$, the
 plan $\l^{s,t}_a$ has $\delta_{g(a)}$ either as first or as second
 marginal, and therefore it is uniquely determined by its marginals
 (it is indeed the product of them). This is the reason why
 condition \eqref{edura} does not show up in the deterministic
 case.

 Finally, we show in Theorem~\ref{optima3} that the two conditions
 are also sufficient, even on general manifolds $D$: if, for some
 $r>1$ and $q\in L^1_{\rm loc}\left((0,T);L^r(D)\right)$, a
 generalized incompressible flow $\eeta$ concentrated on locally
 minimizing curves for the Lagrangian ${\cal L}_q$ satisfies
 $$
 \text{for all $[s,t]\subset (0,T)$, \quad $\l^{s,t}_a$ is
 $c_q^{s,t}$-optimal for $\mu_D$-a.e. $a\in D$,}
 $$
 then $\eeta$ is optimal in $[0,T]$, and $q$ is the pressure field.

 These results show a somehow unexpected connection between the
 variational theory of incompressible flows and the theory developed by
 Bernard-Buffoni \cite{bernard} of measures in the space of action-minimizing
 curves; in this framework one can fit Mather's theory as well as optimal
 transportation problems on manifolds, with a geometric cost. In our case
 the only difference is that the Lagrangian is possibly nonsmoooth (but
 hopefully not so bad), and not given \emph{a priori}, but generated by
 the problem itself. Our approach
 also yields (see Corollary~\ref{vpre}) a new variational characterization
 of the pressure field, as a maximizer of the family of functionals
 (for $[s,t]\subset (0,T)$)
 $$
 q\mapsto\int_{\T^d}W_{c^{s,t}_q}(\eta^s_a,\gamma^t_a)\,d\mu_\T(a),
 \qquad Mq\in L^1\left([s,t]\times\T^d\right),
 $$
 where $\eta^s_a,\,\gamma^t_a$ are the marginals of $\l^{s,t}_a$.

 \smallskip
 \noindent {\bf Acknowledgement.} We warmly thank Yann Brenier for
 the many discussions we had together, and for the constant support
 we had from him. We also thank an anonymous referee,
 who read very carefully this paper, for his comments.

 \section{Notation and preliminary results}

 \noindent {\bf Measure-theoretic notation.} We start by recalling
 some basic facts in Measure Theory. Let $X,\,Y$ be \emph{Polish}
 spaces, i.e. topological spaces whose topology is induced by a
 complete and separable distance. We endow a Polish space $X$ with
 the corresponding Borel $\sigma$-algebra and denote by
 $\Probabilities{X}$ (resp. $\PlusMeasures X$, $\Measures X$) the
 family of Borel probability (resp. nonnegative and finite, real
 and with finite total variation) measures in $X$. For $A\subset X$
 and $\mu\in\Measures{X}$ the \emph{restriction} $\mu\llcorner A$
 of $\mu$ to $A$ is defined by $\mu\llcorner A(B):=\mu(A\cap B)$. We
 will denote by $\ii:X\to X$ the identity map.

 \begin{definition}[Push-forward]
 Let $\mu\in \Measures{X}$ and let $f:X\to Y$ be a Borel map. The
 push-forward $f_\#\mu$ is the measure in $Y$ defined by
 $f_\#\mu(B)=\mu(f^{-1}(B))$ for any Borel set $B\subset Y$. The
 definition obviously extends, componentwise, to vector-valued
 measures.
 \end{definition}

 It is easy to check that $f_\#\mu$ has finite total variation as
 well, and that $|f_\#\mu|\leq f_\#|\mu|$. An elementary
 approximation by simple functions shows the change of variable
 formula
 \begin{equation}\label{chainrule}
 \int_Yg\,df_\#\mu=\int_X g\circ f\,d\mu
 \end{equation}
 for any bounded Borel function (or even either nonnegative or
 nonpositive, and $\overline{\R}$-valued, in the case
 $\mu\in\PlusMeasures X$) $g:Y\to\R$.

 \begin{definition}[Narrow convergence and compactness]\label{dnarrow}
 Narrow (sequential) convergence in $\Probabilities{X}$ is the
 convergence induced by the duality with $C_b(X)$, the space of
 continuous and bounded functions in $X$. By Prokhorov theorem, a
 family $\mathscr F$ in $\Probabilities{X}$ is sequentially
 relatively compact with respect to the narrow convergence if and
 only if it is tight, i.e. for any $\eps>0$ there exists a compact
 set $K\subset X$ such that $\mu(X\setminus K)<\eps$ for any
 $\mu\in \mathscr F$.
 \end{definition}

 In this paper we use only the ``easy'' implication in Prokhorov
 theorem, namely that any tight family is sequentially relatively
 compact. It is immediate to check that a sufficient condition for
 tightness of a family $\mathscr F$ of probability measures is the
 existence of a \emph{coercive} functional $\Psi:X\to [0,+\infty]$
 (i.e. a functional such that its sublevel sets $\{\Psi\leq t\}$,
 $t\in\R^+$, are relatively compact in $X$) such that
 $$
 \int_X\Psi(x)d\mu(x)\leq 1\qquad\forall\mu\in\mathscr F.
 $$
 \begin{lemma}[\cite{amlisa}, Lemma~2.4]\label{ljensen}
 Let $\mu\in\Probabilities{X}$ and $\uu\in L^2(X;\R^m)$. Then, for
 any Borel map $f:X\to Y$, $f_\#(\uu\mu)\ll f_\#\mu$ and its
 density $\vv$ with respect to $f_\#\mu$ satisfies
 $$
 \int_Y|\vv|^2\,df_\#\mu\leq\int_X|\uu|^2\,d\mu.
 $$
 Furthermore, equality holds if and only if $\uu=\vv\circ f$
 $\mu$-a.e. in $X$.
 \end{lemma}

 Given $\mu\in\PlusMeasures{X\times Y}$, we shall denote by
 $\mu_x\otimes\lambda$ its \emph{disintegration} via the projection
 map $\pi(x,y)= x$: here $\lambda=\pi_\#\mu\in\PlusMeasures X$, and
 $x\mapsto\mu_x\in\Probabilities Y$ is a Borel map (i.e.
 $x\mapsto\mu_x(A)$ is Borel for all Borel sets $A\subset Y$)
 characterized, up to $\lambda$-negligible sets, by
 \begin{equation}\label{disinto}
 \int_{X\times Y}f(x,y)\,d\mu(x,y)=\int_X\left(\int_Y
 f(x,y)\,d\mu_x(y)\right)\, d\lambda(x)
 \end{equation}
 for all nonnegative Borel map $f$. Conversely, any $\lambda$ and
 any Borel map $x\mapsto\mu_x\in\Probabilities{Y}$ induce a
 probability measure $\mu$ in $X\times Y$ via \eqref{disinto}.

 \noindent {\bf Function spaces.} We shall denote by $\Omega(D)$
 the space $C([0,T];D)$, and by $\omega:[0,T]\to D$ its typical
 element. The evaluation maps at time $t$,
 $\omega\mapsto\omega(t)$, will be denoted by $e_t$.

 If $D$ is a smooth, compact Riemannian manifold without boundary
 (typically the $d$-dimensional flat torus $\T^d$), we shall denote
 $\mu_D$ its volume measure, and by $d_D$ its Riemannian distance,
 normalizing the Riemannian metric so that $\mu_D$ is a probability
 measure. Although it does not fit exactly in this framework, we
 occasionally consider also the case $D=[0,1]^d$, because many
 results have already been obtained in this particular case.

 We shall often consider measures $\eeta\in\PlusMeasures{\O(D)}$
 such that $(e_t)_\#\eeta\ll\mu_D$; in this case we shall denote by
 $\rho^{\seeta}:[0,T]\times D\to [0,+\infty]$ the density, characterized by
 $$
 \rho^{\seeta}(t,\cdot)\mu_D:=(e_t)_\#\eeta,\qquad t\in [0,T].
 $$

 We denote by $\SDiff{D}$ the measure-preserving diffeomorphisms of
 $D$, and by $S(D)$ the measure-preserving maps in $D$:
 \begin{equation}\label{defSD}
 S(D):=\left\{g:D\to D:\ g_\#\mu_D=\mu_D\right\}.
 \end{equation}
 We also set
 \begin{equation}\label{defSDi}
 S^i(D):=\left\{g\in S(D): \text{$g$ is $\mu_D$-essentially
 injective} \right\}.
 \end{equation}
 For any $g\in S^i(D)$ the inverse $g^{-1}$ is well defined up to
 $\mu_D$-negligible sets, $\mu_D$-measurable, and $g^{-1}\circ
 g=\ii=g\circ g^{-1}$ $\mu_D$-a.e. in $D$. In particular, if $g\in
 S^i(D)$, $g^{-1} \in S^i(D)$.

 We shall also denote by $\Gamma(D)$ the family of
 measure-preserving plans, i.e. the probability measures in
 $D\times D$ whose first and second marginal are $\mu_D$:
 \begin{equation}\label{defGammaD}
 \Gamma(D):=\left\{\gamma\in\Probabilities{D\times D}:\
 (\pi_1)_\#\gamma=\mu_D,\,\, (\pi_2)_\#\gamma=\mu_D\right\}
 \end{equation}
 (here $\pi_1,\,\pi_2$ are the canonical coordinate projections).

 Recall that $\SDiff{D}\subset S^i(D)\subset S(D)$ and that any
 element $g\in S(D)$ canonically induces a measure preserving plan
 $\gamma_g$, defined by
 $$
 \gamma_g:=(\ii\times g)_\#\mu_D.
 $$
 Furthermore, this correspondence is continuous, as long as
 convergence in $L^2(\mu)$ of the maps $g$ and narrow convergence
 of the plans are considered (see for instance Lemma~2.3 in
 \cite{amlisa}). Moreover
 \begin{equation}\label{density1}
 \overline{\left\{\gamma_g:\ g\in S^i(D)\right\}}^{\,{\rm
 narrow}}=\Gamma(D),
 \end{equation}
 \begin{equation}\label{density2}
 \overline{\SDiff{D}}^{L^2(\mu_D)}=S(D) \qquad\text{if $D=[0,1]^d$,
 with $d\geq 2$}
 \end{equation}
 (the first result is standard, see for example the explicit construction in
 \cite[Theorem 1.4 (i)]{brengang} in the case $D=[0,1]^d$, while
 the second one is proved in \cite[Corollary 1.5]{brengang})

 \noindent {\bf The continuity equation.} In the sequel we shall
 often consider weak solutions $\mu_t\in\Probabilities{D}$ of the
 continuity equation
 \begin{equation}\label{conti}
 \p_t \mu_t+\div (\vv_t\mu_t)=0,
 \end{equation}
 where $t\mapsto\mu_t$ is narrowly continuous (this is not
 restrictive, see for instance Lemma~8.1.2 of \cite{amgisa}) and
 $\vv_t(x)$ is a suitable velocity field with
 $\Vert\vv_t\Vert_{L^2(\mu_t)}\in L^1(0,T)$ (formally, $\vv_t$ is a
 section of the tangent bundle and $|\vv_t|$ is computed according
 to the Riemannian metric). The equation is understood in a weak
 (distributional) sense, by requiring that
 $$
 \frac{d}{dt}\int_D\phi(t,x)\,d\mu_t(x)=
 \int_D\p_t\phi+\langle\nabla\phi,\vv_t\rangle\,d\mu_t
 \qquad\text{in ${\mathcal D}'(0,T)$}
 $$
 for any $\phi\in C^1\left((0,T)\times D\right)$
 with bounded first derivatives and support
 contained in $J\times D$, with $J\Subset(0,T)$. In the case when $D\subset\R^d$
 is compact, we shall consider functions
 $\phi\in C^1\left((0,T)\times\R^d\right)$, again
 with support contained in $J\times\R^d$, with $J\Subset(0,T)$.

 The following general principle allows to lift solutions of the
 continuity equation to measures in the space of continuous paths.

 \begin{theorem}[Superposition principle]\label{tsuperp}
 Assume that either $D$ is a compact subset of $\R^d$, or $D$ is a
 smooth compact Riemannian manifold without boundary, and let
 $\mu_t:[0,T]\to\Probabilities{D}$ be a narrowly continuous
 solution of the continuity equation \eqref{conti} for a suitable
 velocity field $\vv(t,x)=\vv_t(x)$ satisfying
 $\Vert\vv_t\Vert^2_{L^2(\mu_t)}\in L^1(0,T)$. Then there exists
 $\eeta\in\Probabilities{\Omega(D)}$ such that
 \begin{itemize}
 \item[(i)] $\mu_t=(e_t)_\#\eeta$ for all $t\in [0,T]$; \item[(ii)]
 the following energy inequality holds:
 $$
 \int_{\Omega(D)}\int_0^T|\dot{\omega}(t)|^2\,dt
 \,d\eeta(\omega)\leq \int_0^T\int_{D}|\vv_t|^2\,d\mu_t\,dt.
 $$
 \end{itemize}
 \end{theorem}

 \begin{Proof} In the case when $D=\R^d$ (and therefore also when
 $D\subset\R^d$ is closed) this result is proved in Theorem~8.2.1
 of \cite{amgisa} (see also \cite{bangert}, \cite{smirnov},
 \cite{bernard1} for related results). In the case when $D$ is a
 smooth, compact Riemannian manifold we recover the same result
 thanks to an isometric embedding in $\R^m$, for $m$ large enough.
 \end{Proof}

 \section{Variational models for generalized geodesics}
 \label{secvarmodel}

 \subsection{Arnold's least action problem}
 Let $f,\,h\in\SDiff{D}$ be given. Following Arnold \cite{arnold},
 we define $\delta^2(f,h)$ by minimizing the action
 $$
 \AA_T(g):=T \int_0^T \int_D \frac{1}{2}|\dot g(t,x)|^2\,d\mu_D(x) \,dt,
 $$
 among all smooth curves
 $$
 [0,T] \ni t \mapsto g(t,\cdot) \in \SDiff{D}
 $$
 connecting $f$ to $h$. By time rescaling, $\delta$ is independent
 of $T$. Since right composition with a given element
 $g\in\SDiff{D}$ does not change the action (as it amounts just to
 a relabelling of the initial position with $g$), the distance
 $\delta$ is right invariant, so it will be often useful to assume,
 in the minimization problem, that $f$ is the identity map.

 The action $\AA_T$ can also be computed in terms of the velocity
 field $\uu$, defined by $\uu(t,x)=\dot g(t,y)|_{y=g^{-1}(t,x)}$, as
 $$
 \AA_T(\uu)=T \int_0^T \int_D \frac{1}{2}|\uu(t,x)|^2\,d\mu_D(x) \,dt.
 $$

 As we mentioned in the introduction, connections between this
 minimization problem and \eqref{Eu} were achieved first by Ebin
 and Marsden, and then by Brenier: in \cite{brenier1},
 \cite{brenier4} he proved that if $(\uu,p)$ is a smooth solution
 of the Euler equation in $[0,T] \times D$, with $D=[0,1]^d$, and
 the inequality in \eqref{goodp} is strict, then the flow $g(t,x)$ of $\uu$ is the unique
 solution of Arnold's minimization problem with $f=\ii$, $h=g(T,\cdot)$.

 By integrating the inequality $d_D^2(h(x),f(x))\leq\int_0^1|\dot
 g(t,x)|^2\,dt$ one immediately obtains that $\Vert
 h-f\Vert_{L^2(D)}\leq\sqrt{2}\delta(f,h)$; Shnirelman proved in
 \cite{shnir2} that in the case $D=[0,1]^d$ with $d\geq 3$ the
 Arnold distance is topologically equivalent to the $L^2$ distance:
 namely, there exist $C>0$, $\a >0$ such that
 \begin{equation}\label{stimaarnold}
 \d(f,g) \leq C \|f-g\|_{L^2(D)}^\a\qquad\forall f,\,g\in\SDiff{D}.
 \end{equation}
 Shnirelman also proved in \cite{shnir1} that when $d\geq 3$ the
 infimum is not attained in general and that, when $d=2$,
 $\delta(\ii,h)$ need not be finite (i.e., there exist
 $h\in\SDiff{D}$ which cannot be connected to $\ii$ by a path with
 finite action).

 \subsection{Brenier's Lagrangian model and its extensions}
 \label{sectBrenierlagr} In \cite{brenier1}, Brenier proposed a
 relaxed version of the Arnold geodesic problem, and here we
 present more general versions of Brenier's relaxed problem,
 allowing first for final data in $\Gamma(D)$, and then for initial
 and final data in $\Gamma(D)$.

 Let $\gamma\in\Gamma(D)$ be given; the class of admissible paths,
 called by Brenier \textit{generalized
 incompressible flows}, is made by the probability measures $\eeta$
 on $\Omega(D)$ such that
 $$
 (e_t)_\# \eeta=\mu_D \qquad \forall t \in [0,T].
 $$
 Then the action of an admissible $\eeta$ is defined as
 $$
 \AA_T(\eeta):=\int_{\O(D)} \AA_T(\omega)\,d\eeta(\omega),
 $$
 where
 \begin{equation}\label{actionomega}
 \AA_T(\omega):=\begin{cases}
 T\int_0^T\frac{1}{2}|\dot\omega(t)|^2\,dt &\text{if $\omega$ is
 absolutely continuous in $[0,T]$}
 \\
 +\infty &\text{otherwise,}
 \end{cases}
 \end{equation}
 and $\overline{\delta}^2(\gamma_\sii,\gamma)$ is defined by
 minimizing $\AA_T(\eeta)$ among all generalized incompressible
 flows $\eeta$ connecting $\gamma_\sii$ to $\gamma$, i.e. those
 satisfying
 \begin{equation}\label{cuoco}
 (e_0,e_T)_\#\eeta=\gamma.
 \end{equation}

 Notice that it is not clear, in this purely Lagrangian
 formulation, how the relaxed distance
 $\overline\delta(\eta,\gamma)$ between two measure preserving
 plans might be defined, not even when $\eta$ and $\gamma$ are
 induced by maps $g,\,h$. Only when $g\in S^i(D)$ we might use the
 right invariance and define
 $\overline{\delta}(\gamma_g,\gamma_h):=
 \overline\delta(\g_\sii,\g_{h\circ
 g^{-1}})$.

 These remarks led us to the following more general problem: let us
 denote
 $$
 \Oa(D):=\O(D)\times D,
 $$
 whose typical element will be denoted by $(\omega,a)$, and let us denote by
 $\pi_D:\Oa(D)\to D$ the canonical projection.
 We consider probability measures $\eeta$ in $\Oa(D)$ having $\mu_D$
 as second marginal, i.e. $(\pi_D)_\#\eeta=\mu_D$; they can be canonically represented as
 $\eeta_a\otimes\mu_D$, where $\eeta_a\in\Probabilities{\O(D)}$. The
 incompressibility constraint now becomes
 \begin{equation}\label{cuoco3}
 \int_D (e_t)_\#\eeta_a\,d\mu_D(a)=\mu_D\qquad\forall t\in [0,T],
 \end{equation}
 or equivalently $(e_t)_\#\eeta=\mu_D$ for all $t$, if we consider
 $e_t$ as a map defined on $\Oa(D)$. Given initial and final data
 $\eta=\eta_a\otimes\mu_D,\,
 \gamma=\gamma_a\otimes\mu_D\in\Gamma(D)$, the constraint
 \eqref{cuoco} now becomes
 \begin{equation}\label{cuoco1}
 (e_0,\pi_D)_\#\eeta=\eta_a\otimes\mu_D,\qquad
 (e_T,\pi_D)_\#\eeta=\gamma_a\otimes\mu_D.
 \end{equation}
 Equivalently, in terms of $\eeta_a$ we
 can write
 \begin{equation}\label{cuoco2}
 (e_0)_\#\eeta_a=\eta_a,\qquad (e_T)_\#\eeta_a=\gamma_a.
 \end{equation}
 Then, we define $\overline\delta^2(\eta,\gamma)$ by minimizing the
 action
 $$
 \int_{\Oa(D)}\AA_T(\omega)\,d\eeta(\omega,a)
 $$
 among all generalized incompressible flows $\eeta$ (according to
 \eqref{cuoco3}) connecting $\eta$ to $\gamma$ (according to
 \eqref{cuoco1} or \eqref{cuoco2}). Notice that $\overline\delta^2$
 is independent of $T$, because the action is scaling invariant; so
 we can use any interval $[a,b]$ in place of $[0,T]$ to define
 $\overline\delta$, and in this case we shall talk of generalized
 flow between $\eta$ and $\gamma$ in $[a,b]$ (this extension will
 play a role in Remark~\ref{rlocalgeo} below).

 When $\eta_a=\delta_a$ (i.e. $\eta=\gamma_{\sii}$), \eqref{cuoco2}
 tells us that almost all trajectories of $\eeta_a$ start from $a$:
 then $\int_D\eeta_a\,d\mu_D(a)$ provides us with a solution of Brenier's
 original model with the same action, connecting $\gamma_\sii$ to
 $\gamma$. Conversely, any solution $\nnu$ of this model can be
 written as $\int_D\nnu_a\,d\mu_D$, with $\nnu_a$ concentrated on
 the curves starting at $a$, and $\nnu_a\otimes\mu_D$
 provides us with an admissible path for our generalized problem,
 connecting $\gamma_\sii$ to $\gamma$, with the same action.

 Let us now analyze the properties of
 $(\Gamma(D),\overline\delta)$; the fact that this is a metric
 space and even a \emph{length space} (i.e. any two points can be
 joined by a geodesic with length equal to the distance) follows by
 the basic operations \emph{reparameterization}, \emph{restriction}
 and \emph{concatenation} of generalized flows, that we are now
 going to describe.

 \begin{remark}[Repameterization]\label{rdavid}
 {\rm Let $\chi:[0,T]\to [0,T]$ be a $C^1$ map with $\dot\chi>0$,
 $\chi(0)=0$ and $\chi(T)=T$. Then, right composition of $\omega$
 with $\chi$ induces a transformation $\eeta\mapsto \chi_\#\eeta$
 between generalized incompressible flows that preserves the
 initial and final conditionl. As a consequence, if $\eeta$ is optimal the functional
 $\chi\mapsto\AA_T(\chi_\#\eeta)$ attains its minimum when
 $\chi(t)=t$. Changing variables we obtain
 $$
 \AA_T(\chi_\#\eeta)=T\int_0^T\dot\chi^2(t)
 \int_{\Oa(D)}\frac{1}{2}|\dot\omega|^2(\chi(t))\,d\eeta(\omega,a)\,dt=
 T\int_0^T\frac{1}{\dot g(s)}
 \int_{\Oa(D)}\frac{1}{2}|\dot\omega|^2(s)\,d\eeta(\omega,a) \,ds
 $$
 with $g=\chi^{-1}$. Therefore, choosing $g(s)=s+\e\phi(s)$, with
 $\phi\in C^1_c(0,T)$, the first variation gives
 $$
 \int_0^T\biggl(
 \int_{\Oa(D)}|\dot\omega|^2(s)\,d\eeta(\omega,a)\biggr)\,\dot\phi(s)
 \,ds=0.
 $$
 This proves that
 $s\mapsto\int_{\Oa(D)}|\dot\omega|^2(s)\,d\eeta(\omega,a)$ is
 equivalent to a constant. We shall call the square root of this quantity
 \emph{speed} of $\eeta$. }\end{remark}

 \begin{remark}[Restriction and concatenation]\label{rlocalgeo}
 {\rm Let $[s,t]\subset [0,T]$ and let $r_{s,t}:C([0,T];D)\to
 C([s,t];D)$ be the restriction map. It is immediate to check that,
 for any generalized incompressible flow $\eeta=\eeta_a\otimes\mu_D$
 in $[0,T]$ between $\eta$ and $\gamma$, the measure
 $(r_{s,t})_\#\eeta$ is a generalized incompressible flow in
 $[s,t]$ between $\eta_s:=(e_s)_\#\eeta_a\otimes\mu_D$ and
 $\gamma_t:=(e_t)_\#\eeta_a\otimes\mu_D$, with action equal to
 $$(t-s)\int_{\Oa(D)}
 \int_s^t\frac{1}{2}|\dot\omega(\tau)|^2\, d\tau\,d\eeta(\omega,a).
 $$
 Let $s<l<t$ and let $\eeta=\mmu_a\otimes\mu_D$,
 $\nnu=\nnu_a\otimes\mu_D$ be generalized incompressible flows,
 respectively defined in $[s,l]$ and $[l,t]$, and joining $\eta$ to
 $\gamma$ and $\gamma$ to $\theta$. Then, writing
 $\gamma_a=(e_l)_\#\eeta_a=(e_l)_\#\nnu_a$, we can disintegrate
 both $\eeta_a$ and $\nnu_a$ with respect to $\gamma_a$ to obtain
 $$
 \eeta_a=\int_D\eeta_{a,x}\,d\gamma_a(x)\in\Probabilities{C([s,l];D)},\qquad
 \nnu_a=\int_D\nnu_{a,x}\,d\gamma_a(x)\in\Probabilities{C([l,t];D)},
 $$
 with $\eeta_{a,x}$, $\nnu_{a,x}$ concentrated on the curves $\omega$ with
 $\omega(l)=x$. We can then consider the image $\llambda_{x,a}$,
 via the concatenation of paths (from the product of $C([s,l];D)$
 and $C([l,t];D)$ to $C([s,t];D)$), of the product measure
 $\eeta_{a,x}\times\nnu_{a,x}$ to obtain a probability measure in
 $C([s,t];D)$ concentrated on paths passing through $x$ at time
 $l$. Eventually, setting
 $$
 \llambda=\int_{D\times D}\llambda_{x,a}\,d(\gamma_a\otimes\mu_D)(x,a),
 $$
 we obtain a generalized incompressible flow in $[s,t]$ joining
 $\eta$ to $\theta$ with action given by
 $$
 \frac{t-s}{l-s}\AA_{[s,l]}(\eeta)+\frac{t-s}{t-l}\AA_{[l,t]}(\nnu),
 $$
 where $\AA_{[s,l]}(\eeta)$ is the action of $\eeta$ in $[s,l]$ and
 $\AA_{[l,t]}(\nnu)$ is the action of $\nnu$ in $[l,t]$ (strictly
 speaking, the action of their restrictions).
 }\end{remark}

 A simple consequence of the previous remarks is that
 $\overline\delta$ is a distance in $\Gamma(D)$ (it suffices to
 concatenate flows with unit speed); in addition, the restriction
 of an optimal incompressible flow $\eeta=\eeta_a\otimes\mu_D$
 between $\eta_a\otimes\mu_D$ and $\gamma_a\otimes\mu_D$ to an
 interval $[s,t]$ is still an optimal incompressible flow in $[s,t]$
 between the plans $(e_s)_\#\eeta_a\otimes\mu_D$ and
 $(e_t)_\#\eeta_a\otimes\mu_D$. This property will be useful in
 Section~\ref{secnec}.

 Another important property of $\overline\delta$ that will be
 useful in Section~\ref{secnec} is its lower semicontinuity with
 respect to the narrow convergence, that we are going to prove in
 the next theorem. Another non-trivial fact is the existence of at
 least one generalized incompressible flow with finite action. In
 \cite[Section 4]{brenier1} Brenier proved the existence of such a
 flow in the case $D=\T^d$. Then in \cite[Section 2]{shnir2}, using
 a (non-injective) Lipschitz measure-preserving map from $\T^d$
 to $[0,1]^d$, Shnirelman produced a flow with finite action also
 in this case (see also \cite[Section 3]{brenier4}). In the next
 theorem we will show how to construct a flow with finite action in
 a compact subset $D$ whenever flows with finite action can be built in $D'$ and
 a possibly non-injective, Lipschitz and measure-preserving map $f:D'\to D$ exists.

 \begin{theorem} \label{bendi}
 Assume that $D\subset\R^d$ is a compact set. Then the infimum in
 the definition of $\overline\delta(\eta,\gamma)$ is achieved,
 \begin{equation}\label{bot}
 \text{$(\eta,\gamma)\mapsto\overline\delta(\eta,\gamma)$ is
 narrowly lower semicontinuous}
 \end{equation}
 and
 \begin{equation}\label{bot2}
 \overline{\delta}(\gamma_\sii,\gamma_h)
 \leq\delta(\ii,h)\qquad\forall h\in\SDiff{D}.
 \end{equation}
 Furthermore, $\sup\limits_{\eta,\,\gamma\in\Gamma(D)}
 \overline\delta(\eta,\gamma)\leq\sqrt{d}$ when either $D=[0,1]^d$
 or $D=\T^d$ and, more generally,
 $$
 \sup_{\gamma\in\Gamma(D)} \overline\delta_D(\gamma_{\sii},\gamma)\leq
 {\rm Lip}(f) \sup_{\gamma'\in\Gamma(D')}
 \overline\delta_{D'}(\gamma_{\sii},\gamma')
 $$
 whenever a Lipschitz measure-preserving map $f:D'\to D$ exists.
 \end{theorem}

 \begin{Proof}
 The inequality $\overline\delta(\gamma_\sii,\gamma_h)
 \leq\delta(\ii,h)$ simply follows by the fact that any smooth flow
 $g$ induces a generalized one, with the same action, by the
 formula $\eeta=\Phi_\#\mu_D$, where $\Phi:D\to\Oa(D)$ is the map
 $x\mapsto (g(\cdot,x),x)$. Assuming that some generalized
 incompressible flow with a finite action between $\eta$ and
 $\gamma$ exists, the existence of an optimal one follows by the
 narrow lower semicontinuity of $\eeta\mapsto\AA_T(\eeta)$ (because
 $\omega\mapsto\AA_T(\omega)$ is lower semicontinuous in $\O(D)$)
 and by the tightness of minimizing sequences (because
 $\AA_T(\omega)$ is coercive in $\O(D)$, by the Ascoli-Arzel\`a
 theorem). A similar argument also proves the lower semicontinuity
 of $(\eta,\gamma)\mapsto\overline\delta(\eta,\gamma)$, as the
 conditions \eqref{cuoco3}, \eqref{cuoco1} are stable under narrow
 convergence (of $\eeta$ and $\eta$, $\gamma$).

 When either $D=[0,1]^d$ or $D=\T^d$, it follows by the explicit
 construction in \cite{brenier1}, \cite{shnir2} that
 $\overline\delta(\gamma_\sii,\gamma_h)\leq\sqrt{d}$ for all $h\in
 S(D)$; by right invariance (see Proposition~\ref{propotopo} below)
 the same estimate holds for $\overline\delta(\gamma_f,\gamma_h)$
 with $f\in S^i(D)$; by density and lower semicontinuity it extends
 to $\overline\delta(\eta,\gamma)$, with
 $\eta,\,\gamma\in\Gamma(D)$.

 Let $f:D'\to D$ be a Lipschitz measure-preserving map and $h\in
 S(D)$; we claim that it suffices to show the existence of $\g' \in
 \G(D')$ such that $(f \times f)_\# \g'=(\ii \times h)_\# \mu_D$.
 Indeed, if this is proved, since $f$ naturally induces by left
 composition a map $F$ from $\Oa(D')$ to $\Oa(D)$ given by
 $(\o(t),a) \mapsto (f(\o(t)),a)$, then to any $\eeta \in \O(D')$
 connecting $\ii$ to $\gamma'$ we can associate $F_\# \eeta$, which
 will be a generalized incompressible flow connecting $\ii$ to $h$.
 By the trivial estimate
 $$
 \AA_T(F_\# \eeta) \leq {\rm Lip}^2(f) \AA_T(\eeta),
 $$
 one obtains $\ov \d_D(\gamma_{\sii},h)\leq
 {\rm Lip}(f)\ov\d_{D'}(\gamma_{\sii},\g')$.
 By density and lower semicontinuity we get
 the estimate on $\overline\delta_D(\gamma_{\sii},\gamma)$ for all
 $\gamma\in\Gamma(D)$.

 Thus, to conclude the proof, we have to construct $\gamma'$. Let us
 consider the disintegration of $\mu_{D'}$ induced by the map $f$,
 that is
 \begin{equation}\label{bella}
 \mu_{D'}=\int_D\mu_y\,d\mu_D(y)
 \end{equation}
 where, for $\mu_D$-a.e. $y$, $\mu_y$ is a probability measure in
 $D'$ concentrated on the compact set $f^{-1}(y)$. We now define
 $\g'$ as
 $$
 \g':=\int_D \mu_y\times\mu_{h(y)}\,d\mu_D(y).
 $$
 Clearly the first marginal of $\g'$ is $\mu_{D'}$; since $h \in
 S(D)$, changing variables in \eqref{bella} one has
 $\mu_{D'}=\int_D \mu_{h(y)}\,d\mu_D(y)$, and so also the second
 marginal of $\g'$ is $\mu_D$. Let us now prove that $(f \times
 f)_\# \g'=(\ii \times h)_\# \mu_D$: for any $\phi \in C_b(D\times
 D)$ we have
 \begin{align*}
 \int_{D \times D} \phi(y,y')\,d(f \times f)_\# \g'(y,y')
 &=\int_{D' \times D'} \phi(f(x),f(x'))\,d\g'(x,x')\\
 &=\int_D \int_{D' \times D'}
 \phi(f(x),f(x'))\,d\mu_y(x)\,d\mu_{h(y)}(x')\,d\mu_D(y)\\
 &=\int_D \phi(y,h(y))\,d\mu_D(y),
 \end{align*}
 where in the last equality we used that $\mu_y$ is concentrated on
 $f^{-1}(y)$ and $\mu_{h(y)}$ is concentrated on $f^{-1}(h(y))$ for
 $\mu_D$-a.e. $y$.
  \end{Proof}

 By \eqref{stimaarnold}, \eqref{bot2} and the narrow lower
 semicontinuity of $\overline\delta(\ii,\cdot)$ we get
 \begin{equation}\label{stimaarnoldbis}
 \overline\delta(\gamma_\sii,h)\leq C\|h-\ii\|_{L^2(D)}^\a\qquad
 \text{if $h\in S(D)$, $D=[0,1]^d$, $d\geq 3$.}
 \end{equation}

 We conclude this section by pointing out some additional
 properties of the metric space $(\Gamma(D),\overline\delta)$.

 \begin{proposition} \label{propotopo}
 $(\Gamma(D),\overline\delta)$ is a complete metric space, whose
 convergence implies narrow convergence. Furthermore, the distance
 $\overline\delta$ is right invariant under the action of $S^i(D)$
 on $\Gamma(D)$. Finally, $\overline\delta$-convergence is strictly
 stronger than narrow convergence and, as a consequence,
 $(\Gamma(D),\overline\delta)$ is not compact.
 \end{proposition}
 \begin{Proof}
 We will prove that $\overline{\delta}(\eta,\gamma)\geq
 W_2(\eta,\gamma)$, where $W_2$ is the quadratic Wasserstein
 distance in $\Probabilities{D\times D}$ (with the quadratic cost
 $c((x_1,x_2),(y_1,y_2))=d_D^2(x_1,y_1)/2+d_D^2(x_2,y_2)/2$);
 as this distance metrizes the narrow convergence,
 this will give the implication between
 $\overline\delta$-convergence and narrow convergence. In order to
 show the inequality $\overline\delta(\eta,\gamma)\geq
 W_2(\eta,\gamma)$ we consider an optimal flow $\eeta_a \otimes \mu_D$
 defined in $[0,1]$; then, denoting by $\omega_a\in\O(D)$ the constant
 path identically equal to $a$, and by
 $\nnu_a\in\Probabilities{C([0,1];D\times D)}$
 the measure $\eeta_a\times\delta_{\omega_a}$,
 the measure $\nnu:=\int_D\nnu_a\,d\mu_D(a)\in\Probabilities{C([0,1];D\times D)}$
 provides a ``dynamical transference plan'' connecting $\eta$ to $\gamma$
 (i.e. $(e_0)_\#\nnu=\eta$, $(e_1)_\#\nnu=\gamma$, see
 \cite[Chapter 7]{villaniSF}) whose action is $\overline\delta^2(\eta,\gamma)$;
 since the action of any dynamical transference plan bounds from above
 $W^2_2(\eta,\gamma)$, the inequality is achieved.

 The completeness of $(\Gamma(D),\overline\delta)$ is a consequence
 of the inequality $\overline\delta\geq W_2$ (so that Cauchy
 sequences in this space are Cauchy sequences for the Wasserstein
 distance), the completeness of the Wasserstein spaces of
 probability measures and the narrow lower semicontinuity of
 $\overline\delta$: we leave the details of the simple proof to the
 reader.

 The right invariance of $\overline\delta$ simply follows by the
 fact that $\eta\circ h=\eta_{h(a)}\otimes\mu_D$, $\gamma\circ
 h=\gamma_{h(a)}\otimes\mu_D$, so that
 $$
 \overline\delta(\eta\circ h,\gamma\circ
 h)\leq\overline\delta(\eta,\gamma),
 $$
 because we can apply the same transformation to any admissible
 flow $\eeta_a \otimes \mu_D$ connecting $\eta$ to $\gamma$,
 producing an admissible flow $\eeta_{h(a)} \otimes \mu_D$ between
 $\eta\circ h$ and $\gamma\circ h$ with the same action. If $h\in
 S^i(D)$ the inequality can be reversed, using $h^{-1}$.

 Now, let us prove the last part of the statement. We first show
 that
 \begin{equation}\label{saddam}
 \frac{1}{2}\int_D d^2_D(f,h)\,d\mu_D
 \leq\overline\delta^2(\gamma_f,\gamma_h)\qquad\forall f,\,h\in
 S(D).
 \end{equation}
 Indeed, considering again an optimal flow $\eeta_a \otimes \mu_D$,
 for $\mu_D$-a.e. $a \in D$ we have
 $$
 \frac{1}{2}d_D^2(f(a),h(a))=W_2^2(\d_{f(a)},\d_{h(a)})\leq
 T\int_{\Omega(D)}\int_0^T\frac{1}{2}|\dot\omega(t)|^2\,dt\,d\eeta_a(\omega),
 $$
 and we need only to integrate this inequality with respect to $a$.
 From \eqref{saddam} we obtain that $S(D)$ is a closed subset of
 $\Gamma(D)$, relative to the distance $\overline\delta$. In
 particular, considering for instance a sequence $(g_n)\subset
 S(D)$ narrowly converging to $\gamma\in\Gamma(D)\setminus S(D)$,
 whose existence is ensured by \eqref{density1}, one proves that
 the two topologies are not equivalent and the space is not
 compact.
 \end{Proof}

 Combining right invariance with \eqref{stimaarnoldbis}, we obtain
 \begin{equation}\label{stimaarnoldter}
 \overline\delta(\gamma_g,\gamma_h)=
 \overline\delta(\gamma_\sii,\gamma_{h\circ g^{-1}})\leq
 C\|g-h\|^\a_{L^2(D)}\qquad\forall h\in S(D),\,\,g\in S^i(D)
 \end{equation}
 if $D=[0,1]^d$ with $d\geq 3$. By the density of $S^i(D)$ in $S(D)$ in the $L^2$ norm
 and the lower semicontinuity of $\ov\d$,
 this inequality still holds when $g\in S(D)$.

 \subsection{Brenier's Eulerian-Lagrangian model}
 \label{secConstr} In \cite{brenier4}, Brenier proposed a second
 possible relaxation of Arnold's problem, motivated by the fact
 that this second relaxation allows for a much more precise
 description of the pressure field, compared to the Lagrangian
 model (see Section~\ref{secnec}).

 Still denoting by $\eta=\eta_a\otimes\mu_D\in\Gamma(D)$,
 $\gamma=\gamma_a\otimes\mu_D\in\Gamma(D)$ the initial and final
 plan, respectively, the idea is to add to the Eulerian variable
 $x$ a Lagrangian one $a$ (which, in the case $\eta=\gamma_{\sii}$,
 simply labels the position of the particle at time $0$) and to
 consider the family of distributional solutions, indexed by $a\in
 D$, of the continuity equation
 \begin{equation}
 \label{transport2model} \p_t c_{t,a} + \div (\vv_{t,a} c_{t,a})=0
 \qquad\text{in ${\mathcal D}'((0,T)\times D)$, \quad for $\mu_D$-a.e.
 $a$,}
 \end{equation}
 with the initial and final conditions
 \begin{equation}\label{terminal}
 c_{0,a}=\eta_a,\qquad c_{T,a}=\gamma_a, \qquad\text{for
 $\mu_D$-a.e. $a$.}
 \end{equation}
 Notice that minimization of the kinetic energy
 $\int_0^T\int_D|\vv_{t,a}|^2\,dc_{t,a}\,dt$ among all possible
 solutions of the continuity equation would give, according to
 \cite{bebr}, the optimal transport problem between $\eta_a$ and
 $\gamma_a$ (for instance, a path of Dirac masses on a geodesic
 connecting $g(a)$ to $h(a)$ if $\eta_a=\d_{g(a)}$,
 $\gamma_a=\d_{h(a)}$). Here, instead, by averaging with respect to
 $a$ we minimize the \emph{mean} kinetic energy
 $$
 \int_D\int_0^T\int_D|\vv_{t,a}|^2\,dc_{t,a}\,dt\,d\mu_D(a)
 $$
 with the only \emph{global} constraint between the family
 $\{c_{t,a}\}$ given by the incompressibility of the flow:
 \begin{equation}
 \label{incompr2model} \int_D c_{t,a}\,d\mu_D(a)=\mu_D\qquad\forall t\in
 [0,T].
 \end{equation}
 It is useful to rewrite this minimization problem in terms of the
 the global measure $c$ in $[0,T]\times D\times D$ and the
 measures $c_t$ in $D\times D$
 $$
 c:=c_{t,a}\otimes (\Leb{1}\times\mu_D),\qquad
 c_t:=c_{t,a}\otimes\mu_D
 $$
 (from whom $c_{t,a}$ can obviously be recovered by
 disintegration), and the velocity field
 $\vv(t,x,a):=\vv_{t,a}(x)$: the action becomes
 $$
 \AA_T(c,\vv):= T \int_0^T \int_{D \times D}
 \frac{1}{2}|\vv(t,x,a)|^2\,dc(t,x,a),
 $$
 while \eqref{transport2model} is easily seen to be equivalent to
 \begin{equation}\label{transport2modelbis}
 \frac{d}{dt}\int_{D\times D} \phi(x,a)\,dc_t(x,a)=\int_{D\times D}
 \langle\nabla_x\phi(x,a),\vv(t,x,a)\rangle\,dc_t(x,a)
 \end{equation}
 for all $\phi\in C_b(D\times D)$ with a bounded gradient with
 respect to the $x$ variable.

 Thus, we can minimize the action on the class of couples
 measures-velocity fields $(c,\vv)$ that satisfy
 (\ref{transport2modelbis}) and (\ref{incompr2model}), with the
 endpoint condition \eqref{terminal}. The existence of a minimum in
 this class can be proved by standard compactness and lower
 semicontinuity arguments (see \cite{brenier4} for details). This
 minimization problem leads to a squared distance between $\eta$
 and $\gamma$, that we shall still denote by
 $\overline\delta^2(\eta,\gamma)$. Our notation is justified by the
 essential equivalence of the two models, proved in the next
 section.

 \section{Equivalence of the two relaxed models}
 \label{secEquivalence}

 In this section we show that the Lagrangian model is equivalent to
 the Eulerian-Lagrangian one, in the sense that minimal values are
 the same, and there is a way (not canonical, in one direction) to
 pass from minimizers of one problem to minimizers of the other
 one.

 \begin{theorem}\label{thmequivalence}
 With the notations of Sections \ref{sectBrenierlagr} and
 \ref{secConstr},
 $$
 \min_{\seeta} \AA_T(\eeta)=\min_{(c,\svv)} \AA_T(c,\vv)
 $$
 for any $\eta,\,\gamma\in\Gamma(D)$. More precisely, any minimizer
 $\eeta$ of the Lagrangian model connecting $\eta$ to $\gamma$
 induces in a canonical way a minimizer $(c,\vv)$ of the
 Eulerian-Lagrangian one, and satisfies for $\Leb{1}$-a.e. $t\in
 [0,T]$ the condition
 \begin{equation}\label{contrained}
 \dot\omega(t)=\vv_{t,a}(e_t(\omega))\qquad\text{for $\eeta$-a.e.
 $(\omega,a)$.}
 \end{equation}
 \end{theorem}

 \begin{Proof} Up to an isometric embedding, we shall assume that
 $D\subset\R^m$ isometrically (this is needed to apply Lemma~\ref{ljensen}).
 If $\eeta=\eeta_a\otimes\mu\in\Probabilities{\Oa(D)}$ is a
 generalized incompressible flow, we denote by $D'\subset D$ a
 Borel set of full measure such that $\AA_T(\eeta_a)<\infty$ for
 all $a\in D'$. For any $a\in D'$ we define
 $$
 c^\seeta_{t,a}:=(e_t)_\# \eeta_a, \qquad \mm^\seeta_{t,a}=(e_t)_\#
 \left(\dot\omega(t)\eeta_a \right).
 $$
 Notice that $\mm^\seeta_{t,a}$ is well defined for $\Leb{1}$-a.e.
 $t$, and absolutely continuous with respect to $c^\seeta_{t,a}$,
 thanks to Lemma~\ref{ljensen}; moreover, denoting by
 $\vv^\seeta_{t,a}$ the density of $\mm^\seeta_{t,a}$ with respect
 to $c^\seeta_{t,a}$, by the same lemma we have
 \begin{equation}\label{prodi}
 \int_D|\vv^\seeta_{t,a}|^2\,dc_{t,a}^\seeta\leq
 \int_{\Omega(D)}|\dot\omega(t)|^2\,d\eeta_a(\omega),
 \end{equation}
 with equality only if $\dot\omega(t)=\vv^{\seeta}_{t,a}(e_t(\omega))$ for
 $\eeta_a$-a.e. $\omega$. Then, we define the global measure and
 velocity by
 $$
 c^\seeta:=c^\seeta_{t,a}\otimes (\Leb{1}\times\mu_D), \qquad
 \vv^\seeta(t,x,a)=\vv^{\seeta}_t(x,a):=\vv^\seeta_{t,a}(x).
 $$
 It is easy to check that $(c^\seeta,\vv^\seeta)$ is admissible:
 indeed, writing $\eta=\eta_a\otimes\mu$,
 $\gamma=\gamma_a\otimes\mu_D$, the conditions
 $(e_0)_\#\eeta_a=\eta_a$ and $(e_T)_\#\eeta_a=\gamma_a$ yield
 $c^\seeta_{0,a}=\eta_a$ and $c^\seeta_{T,a}=\gamma_a$ (for
 $\mu_D$-a.e. $a$).

 This proves that \eqref{terminal} is fulfilled; the
 incompressibility constraint \eqref{incompr2model} simply comes
 from \eqref{cuoco3}. Finally, we check \eqref{transport2model} for
 $a\in D'$; this is equivalent, recalling the definition of
 $\vv_{t,a}$, to
 \begin{equation}
 \frac{d}{dt}\int_{D} \phi(x)\,dc^\seeta_{t,a}(x)=
 \int_D\langle\nabla\phi,\mm^\seeta_{t,a}\rangle,
 \end{equation}
 which in turn corresponds to
 \begin{equation}
 \frac{d}{dt}\int_{\Omega(D)} \phi(\omega(t))\,d\eeta_a(\omega)=
 \int_{\Omega(D)}\langle\nabla\phi(\omega(t)),\dot\omega(t)\rangle\,d\eeta_a(\omega).
 \end{equation}
 This last identity is a direct consequence of an exchange of
 differentiation and integral.

 By integrating \eqref{prodi} in time and with respect to $a$ we
 obtain that $\AA_T(c^\seeta,\vv^\seeta)\leq\AA_T(\eeta)$, and
 equality holds only if \eqref{contrained} holds.

 So, in order to conclude the proof, it remains to find, given a
 couple measure-velocity field $(c,\vv)$ with finite action that
 satisfies (\ref{transport2model}), (\ref{terminal}) and
 (\ref{incompr2model}), an admissible generalized incompressible
 flow $\eeta$ with $\AA_T(\eeta)\leq\AA_T(c,\vv)$. By applying
 Theorem~\ref{tsuperp} to the family of solutions of the continuity
 equations \eqref{transport2model}, we obtain probability measures
 $\eeta_a$ with $(e_t)_\#\eeta_a=c_{t,a}$ and
 \begin{equation}\label{prodi2}
 \int_{\Omega(D)}\int_0^T|\dot\omega(t)|^2\,dt\,d\eeta_a(\omega)\leq
 \int_0^T\int_{D}|\vv(t,x,a)|^2\,dc_{t,a}(x)\,dt.
 \end{equation}
 Then, because of \eqref{incompr2model}, it is easy to check that
 $\eeta:=\eeta_a\otimes\mu_D$ is a generalized incompressible flow,
 and moreover $\eeta$ connects $\eta$ to $\gamma$. By integrating
 \eqref{prodi2} with respect to $a$, we obtain that
 $\AA_T(\eeta)\leq\AA_T(c,\vv)$.
 \end{Proof}

 \section{Comparison of metrics and gap phenomena}\label{sgap}

 Throughout this section we shall assume that $D=[0,1]^d$. In
 \cite{shnir2}, Shnirelman proved when $d\geq 3$ the following
 remarkable approximation theorem for Brenier's generalized
 (Lagrangian) flows:

 \begin{theorem}\label{thmapproxdiffeo}
 If $d\geq 3$, then each generalized incompressible flow $\eeta$
 connecting $\ii$ to $h\in\SDiff{D}$ may be approximated together
 with the action by a sequence of smooth flows $(g_k(t,\cdot))$
 connecting $\ii$ to $h$. More precisely:
 \begin{enumerate}
 \item[(i)] the measures $\eeta_k:=(g_k(\cdot,x))_\#\mu_D$ narrowly
 converge in $\O(D)$ to $\eeta$; \item[(ii)]
 $\AA_T(g_k)=\AA_T(\eeta_k) \to \AA_T(\eeta)$.
 \end{enumerate}
 \end{theorem}

 This result yields, as a byproduct, the identity
 \begin{equation}
 \overline\delta(\gamma_\sii,\gamma_h)=\delta(\ii,h)\qquad\text{
 for all $h\in\SDiff{D}$, $d\geq 3$.}
 \end{equation}
 More generally the relaxed distance $\overline\delta(\eta,\gamma)$
 arising from the Lagrangian model can be compared, at least when
 $\eta=\gamma_\sii$ and the final condition $\gamma$ is induced by
 a map $h\in S(D)$, with the relaxation $\delta_*$ of the Arnold
 distance:
 \begin{equation}\label{delta*}
 \delta_*(h):=\inf\left\{\liminf_{n\to\infty}\delta(\ii,h_n):\
 h_n\in\SDiff{D},\,\,\int_D |h_n-h|^2\,d\mu_D\to 0\right\}.
 \end{equation}
 By \eqref{bot} and \eqref{bot2}, we have
 $\delta_*(h)\geq\overline\delta(\gamma_\sii,\gamma_h)$, and a
 \emph{gap phenomenon} is said to occur if the inequality is
 strict.

 In the case $d=2$, while examples of $h\in\SDiff{D}$ such that
 $\delta(\ii,h)=+\infty$ are known \cite{shnir1}, the nature of
 $\delta_*(h)$ and the possible occurrence of the gap phenomenon
 are not clear.

 In this section we prove the non-occurrence of the gap phenomenon
 when the final condition belongs to $S(D)$, and even when it is a transport
 plan, still under the assumption $d\geq 3$. To this aim, we first
 extend the definition of $\delta_*$ by setting
 \begin{equation}\label{delta*bis}
 \delta_*(\gamma):=\inf\left\{\liminf_{n\to\infty}\delta(\ii,h_n):\
 h_n\in\SDiff{D},\,\,\text{$\gamma_{h_n}\to\gamma$
 narrowly}\right\}.
 \end{equation}
 This extends the previous definition \eqref{delta*}, taking into
 account that $\gamma_{h_n}$ narrowly converge to $\gamma_h$ if and
 only if $h_n\to h$ in $L^2(\mu_D)$ (for instance, this is a simple consequence
 of \cite[Lemma 2.3]{amlisa}).

 \begin{theorem} \label{nogapGamma}
 If $d\geq 3$, then
 $\delta_*(\gamma)=\overline\delta(\gamma_\sii,\gamma)$ for all
 $\gamma\in\Gamma(D)$.
 \end{theorem}

 The proof of the theorem, given at the end of this section, is a
 direct consequence of Theorem~\ref{thmapproxdiffeo} and of the
 following approximation result of generalized incompressible flows
 by measure-preserving maps (possibly not smooth, or not
 injective), valid in \emph{any} number of dimensions.

 \begin{theorem}\label{nogapGamma2}
 Let $\gamma\in\Gamma(D)$. Then, for any probability measure
 $\eeta$ on $\O(D)$ such that
 $$
 (e_t)_\# \eeta=\mu_D\quad \forall t \in [0,T], \qquad \quad
 (e_0,e_T)_\# \eeta= \gamma,
 $$
 and $\AA_T(\eeta)<\infty$, there exists a sequence of flows
 $(g_k(t,\cdot))_{k \in \N}\subset
 W^{1,2}\left([0,T];L^2(D)\right)$ such that:
 \begin{itemize}
 \item[(i)] $g_k(t,\cdot)\in S(D)$ for all $t\in [0,T]$, hence
 $\eeta_k:=(\Phi_{g_k})_\#\mu_D$, with
 $\Phi_{g_k}(x)=g_k(\cdot,x)$, are generalized incompressible
 flows; \item[(ii)] $\eeta_k$ narrowly converge in $\O(D)$ to
 $\eeta$ and $\AA_T(g_k)=\AA_T(\eeta_k) \to \AA_T(\eeta)$.
 \end{itemize}
 \end{theorem}

 \begin{Proof} The first three steps of the proof are more or less the same
 as in the proof of Shnirelman's approximation theorem
 (Theorem~\ref{thmapproxdiffeo} in \cite{shnir2}).

 \textbf{Step 1.} Given $\e>0$ small, consider the affine
 transformation of $D$ into the concentric cube $D_\e$ of size
 $1-4\e$:
 $$
 T_\e(x):=(2\e, \ldots,2\e) + (1-4\e)x.
 $$
 This transformation induces a map $\tilde T_\e$ from $\O(D)$ into
 $C([0,T];D_\e)$ (which is indeed a bijection) given by
 $$
 \tilde T_\e (\o)(t):=T_\e(\o(t)) \qquad \forall \o \in \O(D).
 $$
 Then we define $\tilde\eeta_\e:=(\tilde T_\e)_\# \eeta$, and
 $$
 \eeta_\e:=(1-4\e)^d \tilde\eeta_\e + \eeta_{0,\e},
 $$
 where $\eeta_{0,\e}$ is the ``steady'' flow in $D \setminus D_\e$:
 it consists of all the curves in $D \setminus D_\e$ that do not
 move for $0 \leq t \leq T$. It is then not difficult to prove that
 $\eeta_\e\to\eeta$ narrowly and $\AA_T(\eeta_\e) \to
 \AA_T(\eeta)$, as $\e \to 0$.

 Therefore, by a diagonal argument, it suffices to prove our
 theorem for a measure $\eeta$ which is steady near $\p D$. More
 precisely we can assume that, if $\o(0)$ is in the
 $2\e$-neighborhood of $\p D$, then $\o(t) \equiv \o(0)$ for
 $\eeta$-a.e. $\o$. Moreover, arguing as in Step 1 of the proof of
 the above mentioned approximation theorem in \cite{shnir2}, we can
 assume that the flow does not move for $0 \leq t \leq \e$, that
 is, for $\eeta$-a.e. $\o$, $\o(t) \equiv \o(0)$
 for $0 \leq t \leq \e$.

 \textbf{Step 2.} Let us now consider a family of independent
 random variables $\omega_1,\omega_2,\ldots$ defined in a common
 probability space $(Z,\mathcal Z,P)$, with values in $C([0,T],D)$
 and having the same law $\eeta$. Recall that $\eeta$ is steady
 near $\p D$ and for $0 \leq t \leq \e$, so we can see $\omega_i$
 as random variables with values in the subset of $\O(D)$ given by
 the curves which do not move for $0 \leq t \leq \e$ and in the
 $2\e$-neighbourhood of the $\partial D$. By the law of large
 numbers, the random probability measures in $\O(D)$
 $$
 \nnu_N(z):= \frac{1}{N} \sum_{i=1}^N \d_{\omega_i(z)},\quad\qquad
 z\in Z,
 $$
 narrowly converge to $\eeta$ with probability $1$. Moreover,
 always by the law of large numbers, also
 $$
 \AA_T(\nnu_N(z)) \to \AA_T(\eeta)
 $$
 with probability $1$. Thus, choosing properly $z$, we have
 approximated $\eeta$ with measures $\nnu_N$ concentrated on a
 finite number of trajectories $\omega_i(z)(\cdot)$ which are
 steady in $[0,\e]$ and close to $\partial D$. From now on (as
 typical in Probability theory) the parameter $z$ will be tacitly
 understood.

 \textbf{Step 3.} Let $\varphi \in C_c^\i(\R^d)$ be a smooth radial
 convolution kernel with $\varphi(x)=0$ for $|x| \geq 1$ and
 $\varphi(x)>0$ for $|x|<1$. Given a finite number of trajectories
 $\omega_1, \ldots,\omega_N$ as described is step 2, we define
 $$
 a_i(x):=\frac{1}{\e^d} \varphi\left(
 \frac{x-\omega_i(0)}{\e}\right) \qquad \text{if }
 \dist(\omega_i(0),\p D) \geq \e,
 $$
 $$
 a_i(x):=\frac{1}{\e^d} \sum_{\g \in \G}
 \varphi\left(\frac{x-\g(\omega_i(0))}{\e}\right) \qquad\text{if }
 \dist(\omega_i(0),\p D) \leq \e,
 $$
 where $\G$ is the discrete group of motions in $\R^n$ generated by
 the reflections in the faces of $D$. It is easy to check that
 $\int a_i=1$ and that $\supp(a_i)$ is the intersection of $D$ with
 the closed ball $\overline{B}_\e(\omega_i(0))$. Define
 $$
 g_{i,t}(x):=\omega_i(t) + (x-\omega_i(0)) \qquad \forall
 i=1,\ldots,n.
 $$
 Let
 $\MM_N:=(a_1,\ldots,a_N,g_{1,t}(x),\ldots,g_{N,t}(x))$
 and let us consider the generalized flow $\eeta_N$ associated
 to $\MM_N$, given by
 \begin{equation}
 \label{eqmultiflow} \int_{\O(D)} f(\o)\,d\eeta_{N}:=
 \frac{1}{N}\sum_{i=1}^N\int_{D} a_i(x) f(t \mapsto
 g_{i,t}(x))\,dx
 \end{equation}
 (that is, $\eeta_N$ is the measure in the space of paths given by
 $\frac{1}{N}\sum_i\int_D a_i(x)\d_{g_{i,\cdot}(x)}dx$).
 The measure $\eeta_N$ is well defined for the following reason:
 if $\dist(\omega_i(0),\p D)\leq\e$ we have $g_{i,t}(x)=x$, and
 if $\dist(\omega_i(0),\p D)>\e$ and $a_i(x)>0$ we still have that the curve
 $t\mapsto g_{i,t}(x)$ is contained in $D$
 because $a_i(x)>0$ implies $|x-\omega_i(0)|\leq\e$ and,
 by construction, $\dist(\omega_i(t),\p D)\geq\e$ for all times.
Since the density $\rho^{\seeta_N}$ induced by $\eeta_N$ is given
by
 $$
 \rho^N(t,x):=\frac{1}{N}\sum_{i=1}^N a_i(x+\omega_i(0)-\omega_i(t)),
 $$
 the flow $\eeta_N$ is not measure preserving.
 However we are more or less in the same situation as in Step 3 in the proof of
 the approximation theorem in \cite{shnir2} (the only difference
 being that we do not impose any final data). Thus, by \cite[Lemma
 1.2]{shnir2}, with probability $1$
 \begin{equation}
 \label{eqLargeNumberLaw}
 \begin{split}
 &\sup_{x,t} |\rho^N(t,x) -1| \to 0, \\
 &\sup_{x,t} |\p^\a_x\rho^N(t,x)| \to 0 \quad \forall \a, \\
 &\int_D \int_0^T |\p_t\rho^N(t,x)|^2 \,dt\,dx \to 0
 \end{split}
 \end{equation}
 as $N \to \infty$. By the first two equations in
 (\ref{eqLargeNumberLaw}), we can left compose $g_{i,t}$ with a
 smooth correcting flow $\zeta_t^N(x)$ as in Step 3 in the proof of
 the approximation theorem in \cite{shnir2}, in such a way that the
 flow $\tilde\eeta_N$ associated to $\tilde\MM_N:=(a_1,\ldots,a_N,\zeta_t^N \circ
 g_{1,t}(x),\ldots,\zeta_t^N \circ g_{N,t}(x))$ via the formula
 analogous to (\ref{eqmultiflow}) is incompressible. Moreover, thanks to the
 third equation in \eqref{eqLargeNumberLaw} and the convergence of
 $\AA_T(\nnu_N)$ to $\AA_T(\eeta)$, one can prove that
 $\AA_T(\tilde\eeta_N)\to\AA_T(\eeta)$ with probability $1$.

 We observe that, since $\eeta$ is steady for $0 \leq t \leq \e$,
 the same holds by construction for $\tilde\eeta_N$. Without
 loss of generality, we can therefore assume that $\zeta_t^N$ does
 not depend on $t$ for $t \in [0,\e]$.

 \textbf{Step 4.} In order to conclude, we see that the only
 problem now is that the flow $\tilde\eeta_N$ associated to
 $\tilde\MM_N$ is still
 non-deterministic, since if $x \in \supp(a_i) \cap \supp(a_j)$ for
 $i \neq j$, then more that one curve starts from $x$. Let us
 partition $D$ in the following way:
 $$
 D=D_1 \cup D_2 \cup \ldots \cup D_L\cup E,
 $$
 where $E$ is $\Leb{d}$-negligible, any set $D_j$ is open, and all $x\in
 D_j$ belong to the interior of the supports of exactly $M=M(j)\leq N$
 sets $a_i$, indexed by $1\leq i_1<\cdots<i_M\leq N$ (therefore $L \leq 2^N$). This
 decomposition is possible, as $E$ is contained in the union of the
 boundaries of $\supp a_i$, which is $\Leb{d}$-negligible.

 Fix one of the sets $D_j$ and assume just for notational
 simplicity that $i_k=k$ for $1\leq k\leq M$. We are going to
 modify the flow $\tilde\eeta_N$ in $D_j$, increasing a
 little bit its action (say, by an amount $\alpha>0$), in such a
 way that for each point in $D_j$ only one curve starts from it.
 Given $x\in D_j$, we know that
 $M$ curves start from it, weighted with mass $a_k(x)>0$, and
 $\sum_{k=1}^M a_k(x)=1$. These curves coincide for $0 \leq t \leq \e$
 (since nothing moves), and then separate. We want to partition
 $D_j$ in $M$ sets $E_k$, with
 $$
 \Leb{d}(E_k)=\int_{D_j} a_k(x)\,dx,\quad \qquad 1\leq k\leq M
 $$
 in such a way that, for any $x \in E_k$, only one curve $\o^k_x$
 starts from it at time $0$, $\o^k_x(t)\in D_j$ for $0\leq
 t\leq\e$, and the map $E_k \ni x \mapsto \o^k_x(\e)\in D_j$ pushes forward
 $\Leb{d}\llcorner E_k$ into $a_k\Leb{d}\llcorner D_j$.
 Moreover, we want the incompressibility condition to be preserved
 for all $t \in [0,\e]$. If this is possible, the proof will be
 concluded by gluing $\o^k_x$ with the only curve starting from
 $\omega^k_x(\e)$ with weight $a_k(\omega^k(\e))$.

 The above construction can be achieved in the following way. First we write the
 interior of $D_j$, up to null measure sets, as a countable union
 of disjoints open cubes $(C_i)$ with size $\delta_i$ satisfying
 \begin{equation}\label{fine}
 \frac{M^2}{\e}\sum_i\frac{\delta_i^2}{\bar
 b_i^2}\Leb{d}(C_i)\leq\alpha,
 \end{equation}
 with $\bar b_i:=\min\limits_{1\leq k\leq M}\min\limits_{C_i}a_k$.
 This is done just considering the union of the grids in $\R^d$
 given by $\Z^d/2^n$ for $n \in \N$, and taking initially our cubes
 in this family; if \eqref{fine} does not hold, we keep splitting
 the cubes until it is satisfied ($\bar b_i$ can only increase
 under this additional splitting, therefore a factor 4 is gained in
 each splitting). Once this partition is given, the idea is to move
 the mass within each $C_i$ for $0 \leq t \leq \e$. At least
 heuristically, one can imagine that in $C_i$ the functions $a_k$
 are almost constant and that the velocity of a generic path in
 $C_i$ is at most of order $\d_i/\e$. Thus, the total energy
 of the new incompressible fluid in the interval $[0,\e]$ will be
 of order
 $$
 \sum_i\int_{C_i} \int_0^\e |\dot \o_x(t)|^2\,dt\,dx \leq
 \frac{C}{\e} \sum_i\d_i^2\Leb{d}(C_i)
 $$
 and the conclusion will follow by our choice of $\delta_i$.

 So, in order to make this argument rigorous, let us fix $i$ and
 let us see how to construct our modified flow in the cube $C_i$
 for $t \in [0,\e]$. Slicing $C_i$ with respect to the first
 $(d-1)$-variables, we see that the transport problem can be solved
 in each slice. Specifically, if $C_i$ is of the form $x^i +
 (0,\delta_i)^d$, and we define
 $$
 m^k:=\int_{C_i} a_k(x)\,dx, \qquad k=1,\ldots,M,
 $$
 whose sum is $\delta_i^d$, then the points which belong to
 $C_i^k:=x^i + (0,\delta_i)^{d-1} \times J_k$ have to move along
 curves in order to push forward $\Leb{d}\llcorner C_i^k$ into
 $a_k\Leb{d}\llcorner C_i$, where $J_k$ are $M$ consecutive open
 intervals in $(0,\delta_i)$ with length $\delta_i^{1-d}m^k$.
 Moreover, this has to be done preserving the incompressibility
 condition.

 If we write $x=(x',x_d) \in \R^d$ with $x'=(x_1,\ldots,x_{d-1})$,
 we can transport the $M$ uniform densities
 $$
 \Haus{1}\llcorner \bigl(x^i + \{x'\}\times J_k \bigr) \quad
 \text{with } x' \in [0,\delta_i]^{(d-1)},
 $$
 into the $M$ densities
 $$
 a_k(x',\cdot)\Haus{1}\llcorner \bigl(x^i + \{x'\}\times
 [0,\delta_i]\bigr)
 $$
 moving the curves only in the $d$-th direction, i.e. keeping $x'$
 fixed. Thanks to Lemma~\ref{chefatica} below and a scaling
 argument, we can do this construction paying at most $M^2\bar
 b_i^{-2}\d_i^3/\e$ in each slice of $C_i$, and therefore with
 a total cost less than
 $$
 \frac{M^2}{\e}\sum_i\frac{\delta_i^{d+2}}{\bar b_i^2}\leq
 \alpha.
 $$
 This concludes our construction.
 \end{Proof}

 \begin{lemma} \label{chefatica} Let $M\geq 1$ be an integer and let
 $b_1,\ldots,b_M:[0,1]\to (0,1]$ be continuous with $\sum_1^M b_k=1$. Setting
 $l_k=\int_0^1b_k\,dt\in (0,1]$, and denoting by $J_1,\ldots,J_M$
 consecutive intervals of $(0,1)$ with length $l_k$, there exists a
 family of uniformly Lipschitz maps $h(\cdot,x)$, with
 $h(t,\cdot)\in S([0,1])$, such that
 $$
 h(1,\cdot)_\#(\chi_{J_k}\Leb{1})=b_k\Leb{1},\qquad k=1,\ldots,M
 $$
 and
 \begin{equation}\label{acest}
 \AA_1(h)\leq \frac{M^2}{\bar{b}^2} ,\qquad\text{
 with $\bar b:=\min\limits_{1\leq k\leq M}\min\limits_{[0,1]} b_k>0$.}
 \end{equation}
 \end{lemma}
 \begin{Proof} We start with a preliminary remark:
 let $J\subset (0,1)$ be an interval with length $l$ and assume
 that $t \mapsto \rho_t$ is a nonnegative Lipschitz map between $[0,1]$ and
 $L^1(0,1)$, with $\rho_t\leq 1$ and $\int_0^1\rho_t\,dx=l$ for all $t \in [0,1]$,
 and let $f(t,\cdot)$ be the unique (on $J$, up to countable sets)
 nondecreasing map pushing $\chi_{J}\Leb{1}$ to $\rho_t$. Assume
 also that $\supp\rho_t$ is an interval and $\rho_t\geq r$
 $\Leb{1}$-a.e. on $\supp\rho_t$, with $r>0$. Under this extra
 assumption, $f(t,x)$ is uniquely determined for all $x\in J$, and
 implicitly characterized by the conditions
 $$
 \int_0^{f(t,x)} \rho_t(y)\,dy=\Leb{1}((0,x)\cap J), \qquad
 f(t,x)\in\supp\rho_t.
 $$
 This implies, in particular, that $f(\cdot,x)$ is continuous
 \emph{for all} $x\in J$. We are going to prove that this map is
 even Lipschitz continuous in $[0,1]$ and
 \begin{equation}\label{stimave}
 |\frac{d}{dt}f(t,x)|\leq\frac{{\rm Lip}(\rho_{\cdot})}{r}
 \qquad\text{for $\Leb{1}$-a.e. $t\in [0,1]$}
 \end{equation}
 for all $x\in J$. To prove this fact, we first notice that the
 endpoints of the interval $\supp\rho_t$ (whose length is at least
 $l$) move at most with velocity ${\rm Lip}(\rho_{\cdot})/r$; then,
 we fix $x\in J=[a,b]$ and consider separately the cases
 $$
 x \in \p J=\{a,b\}, \qquad x \in {\rm Int}(J)=(a,b).
 $$
 In the first case, since for any $t \in [0,1]$
 $$
 \int_0^{f(t,a)} \rho_t(y)\,dy=0, \qquad \int_0^{f(t,b)}
 \rho_t(y)\,dy=\Leb{1}(J),
 $$
 and by assumption $f(t,x)\in\supp\rho_t$ for any $x \in J$, we get
 $\supp\rho_t=[f(t,a),f(t,b)]$ for all $t \in [0,1]$. This,
 together with the fact that the endpoints of the interval
 $\supp\rho_t$ move at most with velocity ${\rm
 Lip}(\rho_{\cdot})/r$, implies (\ref{stimave}) if $x \in \p J$. In
 the second case we have
 $$
 \int_0^{f(t,x)} \rho_t(y)\,dy \in (0,\Leb{1}(J)),
 $$
 therefore $f(t,x)\in {\rm Int}(\supp\rho_t)$ for all $t \in
 [0,1]$. It suffices now to find a Lipschitz estimate of
 $|f(s,x)-f(t,x)|$ when $s,\,t$ are sufficiently close. Assume that
 $f(s,x)\leq f(t,x)$: adding and subtracting
 $\int_0^{f(s,x)}\rho_t(y)\,dy$ in the identity
 $$
 \int_0^{f(t,x)} \rho_t(y)\,dy=\int_0^{f(s,x)}\rho_s(y)\,dy
 $$
 we obtain
 $$
 \int_{f(s,x)}^{f(t,x)}\rho_t(y)\,dy=
 \int_0^{f(s,x)}\rho_s(y)-\rho_t(y)\,dy.
 $$
 Now, as $f(s,x)$ belongs to $\supp\rho_t$ for $|s-t|$ sufficiently
 small, we get
 $$
 r|f(s,x)-f(t,x)|\leq {\rm Lip}(\rho_{\cdot})|t-s|.
 $$
 This proves the Lipschitz continuity of $f(\cdot,x)$ and
 \eqref{stimave}.

 Given this observation, to prove the lemma it suffices to find
 maps $t \mapsto \rho^k_t$ connecting $\chi_{J_k}\Leb{1}$ to $b_k\Leb{1}$
 satisfying:
 \begin{itemize}
 \item[(i)] $\supp\rho^k_t$ is an interval, and $\rho^k_t\geq\min\limits_{[0,1]}
 b_k \geq \bar b$ $\Leb{1}$-a.e. on its support;
 \item[(ii)]
 ${\rm Lip}(\rho^k_{\cdot})\leq \frac{M-1}{2}$ on $[0,\frac{1}{2}]$, and ${\rm Lip}(\rho^k_{\cdot})\leq 2$ on $[\frac{1}{2},1]$;
 \item[(iii)]
 $\sum\limits_{k=1}^M\rho^k_t=1$ for all $t\in [0,1]$.
 \end{itemize}
 Indeed, this would produce maps with time derivative bounded by
 $(M-1)/(2\bar{b})$ on $[0,\frac{1}{2}]$ and bounded by $2/\bar{b}$ on $[\frac{1}{2},1]$,
 and this easily gives \eqref{acest}.

 The construction can be achieved in two steps. First, we connect
 $\chi_{J_k}\Leb{1}$ to $l_k\Leb{1}$ in the time interval
 $[0,\frac{1}{2}]$; then, we connect $l_k\Leb{1}$ to $b_k\Leb{1}$ in
 $[\frac{1}{2},1]$ by a linear interpolation. The Lipschitz constants of
 the second step are easily seen to be less than 2, so let us focus
 on the first interpolation.

 Let us first consider the case of two densities $\rho^1=\chi_{J_1}$ and $\rho^2=\chi_{J_2}$, with $J_1=(0,l_1)$ and $J_2=(l_1,l)$.
 In the time interval $[0,\tau]$, we define the expanding intervals
 $$
 J_{1,t}=(0,l_1 + \frac{t}{\tau}l_2),\quad J_{2,t}=(l_1 - \frac{t}{\tau}l_1,1),
 $$
 so that $J_{k,\tau}=(0,l)$ for $k=1,\,2$, and then define
 $$
 \rho^1_t:=
 \begin{cases}
 1&\text{on $(0,l_1-\frac{t}{\tau}l_1)$,}
 \\
 l_1/l &\text{on $(l_1-\frac{t}{\tau}l_1,l_1 + \frac{t}{\tau}l_2)$,}
 \\
 0 &\text{otherwise}.
 \end{cases}
 \qquad
 \rho^2_t:=
 \begin{cases}
 1&\text{on $(l_1+\frac{t}{\tau}l_2,l)$,}
 \\
 l_2/l &\text{on $(l_1-\frac{t}{\tau}l_1,l_1 + \frac{t}{\tau}l_2)$,}
 \\
 0 &\text{otherwise}.
 \end{cases}
 $$
 By construction $\rho^k_t\geq l_k$ on $J_{k,t}$ for $k=1,\,2$, $\rho^1_t+\rho^2_t=1$,
 and it is easy to see that
 \begin{equation}
 \label{eqboundLip}
 {\rm Lip}(\rho^k_{\cdot})\leq \frac{l_1l_2}{\tau l}\leq \frac{l}{4\tau}.
 \end{equation}

 We can now define the desired interpolation on $[0,\frac{1}{2}]$ for general $M \geq 2$.
 Let us define
 $$
 t_i:=\frac{i}{2(M-1)} \qquad \text{for }i=1,\ldots,M-1,
 $$
 so that $t_{M-1}=\frac{1}{2}$.
 We will achieve our construction of $\rho_t^k$ on $[0,\frac{1}{2}]$ in $M-1$
 steps, where at each step we will progressively define $\rho_t^k$ on
 the time interval $[t_{i-1},t_i]$.

 First, in the time interval $[0,t_1]$,
 we leave fixed $\rho^k_0:=\chi_{J_k}\Leb{1}$ for $k \geq 3$ (if such $k$ exist),
 while we apply the above construction in $J_1 \cup J_2$ to $\rho^1$ and $\rho^2$.
 In this way, on $[0,t_1]$, $\rho^1_0:=\chi_{J_1}\Leb{1}$ is connected to $\rho^1_{t_1}:=\frac{l_1}{l_1+l_2}\chi_{J_1 \cup J_2}\Leb{1}$,
 and $\rho^2_0:=\chi_{J_2}\Leb{1}$ is connected to $\rho^2_{t_1}:=\frac{l_2}{l_1+l_2}\chi_{J_1 \cup J_2}\Leb{1}$.

 Now, as a second step, we want to connect $\rho^k_{t_1}$ to $\frac{l_k}{l_1+l_2+l_3}\chi_{J_1 \cup J_2 \cup J_3}\Leb{1}$
 for $k=1,\,2,\,3$, leaving the other densities fixed. To this aim, we define $\rho^{12}_{t_1}:=\rho^1_{t_1}+\rho^2_{t_1}=\chi_{J_1 \cup J_2}\Leb{1}$.
 In the time interval $[t_1,t_2]$, we leave fixed $\rho^k_0:=\chi_{J_k}\Leb{1}$ for $k \geq 4$ (if such $k$ exist),
 and we apply again the above construction in $J_1 \cup J_2 \cup J_3$ to $\rho^{12}_{t_1}$ and $\rho^3_{t_1}=\chi_{J_3}\Leb{1}$.
 In this way, on $[t_1,t_2]$, $\rho^{12}_{t_1}$ is connected to
 $\rho^{12}_{t_2}:=\frac{l_1+l_2}{l_1+l_2+l_3}\chi_{J_1 \cup J_2 \cup J_3}\Leb{1}$,
 and $\rho^3_{t_1}$ is connected to $\rho^3_{t_2}:=\frac{l_3}{l_1+l_2+l_3}\chi_{J_1 \cup J_2 \cup J_3}\Leb{1}$.
 Finally, it suffices to define $\rho^1_t:=\frac{l_1}{l_1+l_2}\rho^{12}_{t}$ and $\rho^2_t:=\frac{l_2}{l_1+l_2}\rho^{12}_{t}$.

 In the third step we leave fixed the densities $\rho^k_{t_2}$
 for $k \geq 5$, and we do the same construction as before adding the
 first three densities (that is, in this case one defines $\rho^{123}_{t_2}:=\rho^1_{t_2}+\rho^2_{t_2}+\rho^3_{t_2}=\chi_{J_1 \cup J_2 \cup J_3}\Leb{1}$).
 In this way, we connect $\rho^{123}_{t_2}$ to
 $\rho^{123}_{t_3}:=\frac{l_1+l_2+l_3}{l_1+l_2+l_3+l_4}\chi_{J_1 \cup J_2 \cup J_3 \cup J_4}\Leb{1}$
 and $\rho^4_{t_2}$ to $\rho^4_{t_3}:=\frac{l_4}{l_1+l_2+l_3+l_4}\chi_{J_1 \cup J_2 \cup J_3 \cup
 J_4}\Leb{1}$, and then we define
 $\rho^k_t:=\frac{l_k}{l_1+l_2+l_3}\rho^{123}_{t}$ for $k=1,\,2,\,3$.

 Iterating this construction on $[t_i,t_{i+1}]$ for $i \geq 4$, one obtains the desired maps $t\mapsto \rho_t^k$.
 Indeed, by construction $\rho^k_t\geq l_k$ on $J_{k,t}$, and $\sum_{k=1}^M \rho^k_t=1$.
 Moreover, by (\ref{eqboundLip}), it is simple to see that in each time interval $[t_i,t_{i+1}]$ one has the bound
 $$
 {\rm Lip}(\rho^k_{\cdot})\leq \frac{M-1}{2}.
 $$
 So the energy can be easily bounded by
 $1/\bar{b}^2\bigl(\frac{(M-1)^2}{16} + 1\bigr) \leq M^2/\bar{b}^2$.
 \end{Proof}

 \begin{Proof} (of Theorem~\ref{nogapGamma}) By applying Theorem~\ref{nogapGamma2}
 to the optimal $\eeta$ connecting $\ii$ to $\gamma$, we can find
 maps $g_k\in S(D)$ such that $\gamma_{g_k}\to\gamma$ narrowly and
 $$
 \limsup_{k\to\infty}\overline\delta(\gamma_{\sii},\gamma_{g_k})
 \leq\overline\delta(\gamma_{\sii},\gamma).
 $$
 Now, if $d\geq 3$ we can use \eqref{stimaarnoldter}, the triangle
 inequality, and the density of $\SDiff{D}$ in $S(D)$ in the $L^2$ norm, to find maps
 $h_k\in\SDiff{D}$ such that
 $$
 \limsup_{k\to\infty}\overline\delta(\gamma_{\sii},\gamma_{h_k})
 \leq\overline\delta(\gamma_\sii,\gamma)
 $$
 and $\gamma_{h_k}\to\gamma$ narrowly. This gives the thesis.
 \end{Proof}

 \section{Necessary and sufficient optimality conditions}
 \label{secnec}

 In this section we study necessary and sufficient optimality
 conditions for the generalized geodesics; we shall work mainly
 with the Lagrangian model, but we will use the equivalent
 Eulerian-Lagrangian model to transfer regularity informations for
 the pressure field to the
 Lagrangian model. Without any loss of generality, we assume
 throughout this section that $T=1$.

 The pressure field $p$ can be identified, at least as a
 distribution (precisely, an element of the dual of
 $C^1\left([0,1]\times D\right)$), by the so-called dual least
 action principle introduced in \cite{brenier2}. In order to
 describe it, let us build a natural class of first variations
 in the Lagrangian model: given a smooth vector field $\ww(t,x)$,
 vanishing for $t$ sufficiently close to $0$ and $1$, we may define
 the maps $S^\e:\Oa(D)\to\Oa(D)$ by
 \begin{equation}\label{variazw}
 S^\e(\omega,a)(t):=\left(e^{\e\sww_t}\omega(t),a\right),
 \end{equation}
 where $e^{\e\sww_t}x$ is the flow, in the $(\e,x)$ variables,
 generated by the autonomous field $\ww_t(x)=\ww(t,x)$ (i.e.
 $e^{0\sww_t}=\ii$ and
 $\frac{d}{d\e}e^{\e\sww_t}x=\ww(t,e^{\e\sww_t}x)$), and the
 perturbed generalized flows $\eeta_\e:=(S^\e)_\#\eeta$. Notice
 that $\eeta_\e$ is incompressible if $\div\,\ww_t=0$, and more
 generally the density $\rho^{\seeta_\e}$ satisfies for all times
 $t\in (0,1)$ the continuity equation
 \begin{equation}\label{contine}
 \frac{d}{d\e}\rho^{\seeta_\e}(t,x)+\div
 (\ww_t(x)\rho^{\seeta_\e}(t,x))=0.
 \end{equation}
 This motivates the following definition.

 \begin{definition}[Almost incompressible flows]
 \label{defslighcomprflows} We say that a probability measure
 $\nnu$ on $\O(D)$ is a \textit{almost incompressible generalized
 flow} if $\rho^\snnu\in C^1\left([0,1]\times D\right)$ and
 $$
 \Vert\rho^\snnu-1\Vert_{C^1([0,1]\times D)}\leq\frac{1}{2}.
 $$
 \end{definition}

 Now we provide a slightly simpler proof of the characterization
 given in \cite{brenier2} of the pressure field (the original proof
 therein involved a time discretization argument).

 \begin{theorem}
 \label{thmslighcomprflows} For all $\eta,\,\gamma\in\Gamma(D)$
 there exists $p\in \left[C^1([0,1]\times D)\right]^*$ such that
 \begin{equation}\label{augmented}
 \langle p,\rho^\snnu - 1 \rangle_{(C^1)^*,C^1} \leq \AA_1(\nnu)
 - \overline\delta^2(\eta,\gamma)
 \end{equation}
 for all almost incompressible flows $\nnu$ satisfying
 \eqref{cuoco1}.
 \end{theorem}

 \begin{Proof}
 Let us define the closed convex set $C:=\{\rho \in C^1([0,1]\times D):\
 \|\rho-1\|_{C^1}\leq \frac{1}{2}\}$, and the function
 $\phi:C^1\left([0,1]\times D\right) \to \R^+ \cup\{+\i \}$ given by
 $$
 \phi(\rho):=\left\{
 \begin{array}{ll}
 \inf\left\{\AA_1(\nnu):\ \text{$\rho^\snnu=\rho$ and
 \eqref{cuoco1} holds}\right\} &
 \text{if } \rho \in C;\\
 +\infty & \text{otherwise}.
 \end{array}
 \right.
 $$
 We observe that $\phi(1)=\overline\delta^2(\eta,\gamma)$.
 Moreover, it is a simple exercise to prove that $\phi$ is convex
 and lower semicontinuous in $C^1\left([0,1]\times D\right)$. Let
 us now prove that $\phi$ has bounded (descending) slope at $1$,
 i.e.
 $$
 \limsup_{\rho\to 1}\frac{[\phi(1)-\phi(\rho)]^+} {\lVert
 1-\rho\rVert_{C^1}}<+\infty,
 $$
 By \cite[Proposition 2.1]{brenier2} we
 know that there exist $0 <\e<\frac{1}{2}$ and $c>0$ such that, for
 any $\rho\in C$ with $\lVert\rho-1\rVert_{C^1}\leq\e$, there is a
 Lipschitz family of diffeomorphisms $g_\rho(t,\cdot):D \to D$ such
 that
 $$
 g_\rho(t,\cdot)_\#\mu_D=\rho(t,\cdot)\mu_D,
 $$
 $g_\rho(t,\cdot)=\ii$ for $t=0,\,1$,
 and the Lipschitz constant of $(t,x) \mapsto g_\rho(t,x) - x$
 is bounded by $c$. Thus, adapting the
 construction in \cite[Proposition 2.1]{brenier2} (made for
 probability measures in $\O(D)$, and not in $\Oa(D)$), for
 any incompressible flow $\eeta$ connecting $\eta$ to $\gamma$, and any
 $\rho \in C$, we can define an almost incompressible flow $\nnu$
 still connecting $\eta$ to $\gamma$ such that $\rho^\snnu=\rho$,
 and
 $$
 \AA_1(\nnu) \leq \AA_1(\eeta)+c'\lVert \rho-1
 \rVert_{C^1}(1+\AA_1(\eeta)),
 $$
 where $c'$ depends only on $c$ (for instance, we define
 $\nnu:=G_\#\eeta$, where $G:\Oa(D) \to \Oa(D)$ is the map
 induced by $g_\rho$ via the formula $(\o(t),a) \mapsto
 (g_\rho(t,\o(t)),a)$). In particular, considering an optimal $\eeta$,
 we get
 \begin{equation}
 \label{slopephi} \phi(\rho) \leq \phi(1) + c\lVert \rho-1
 \rVert_{C^1}(1+\overline\delta^2(\eta,\gamma))
 \end{equation}
 for any $\rho\in C$ with $\lVert \rho-1 \rVert_{C^1}\leq\e$. This fact
 implies that $\phi$ is bounded on a neighbourhood of $1$ in $C$.
 Now, it is a standard fact of convex analysis that a convex function
 bounded on a convex set is locally Lipschitz on that set.
 This provides the bounded slope property. By a simple
 application of the Hahn-Banach theorem (see for instance
 Proposition~1.4.4 in \cite{amgisa}), it follows that the
 subdifferential of $\phi$ at $1$ is not empty, that is, there
 exists $p$ in the dual of $C^1$ such that
 $$
 \langle p,\rho-1\rangle_{(C^1)^*,C^1}\leq\phi(\rho)-\phi(1).
 $$
 This is indeed equivalent to \eqref{augmented}.
 \end{Proof}

 This result tells us that, if $\eeta$ is an optimal incompressible
 generalized flow connecting $\eta$ to $\gamma$ (i.e.
 $\AA_1(\eeta)=\overline\delta^2(\eta,\gamma)$), and if we consider
 the augmented action
 \begin{equation}\label{AApT}
 \AA^p_1(\nnu):= \int_{\Oa(D)} \int_0^1
 \frac{1}{2}|\dot{\omega}(t)|^2\,dt \,d\nnu(\omega,a) - \langle
 p,\rho_\snnu-1\rangle,
 \end{equation}
 then $\eeta$ minimizes the new action among all almost
 incompressible flows $\nnu$ between $\eta$ and $\gamma$.

 Then, using the identities
 $$
 \frac{d}{d\e} \frac{d}{dt} S^\e(\omega)(t)\biggl\vert_{\e=0}=
 \frac{d}{dt}\ww(t,\omega(t))=\partial_t\ww(t,\omega(t))+
 \nabla_x\ww(t,\omega(t))\cdot\dot\omega(t)
 $$
 and the convergence in the sense of distributions (ensured by
 \eqref{contine}) of $(\rho^{\seeta_\e}-1)/\e$ to $-\div\,\ww$ as
 $\e\downarrow 0$, we obtain
 \begin{equation}\label{befana}
 0=\frac{d}{d\e}\AA_1^p(\eeta_\e)\biggl\vert_{\e=0}=
 \int_{\Oa(D)}\int_0^1\dot\omega(t)\cdot\frac{d}{dt}
 \ww(t,\omega(t))\,dt\,d\eeta(\omega,a)+ \langle
 p,\div\,\ww\rangle.
 \end{equation}
 As noticed in \cite{brenier2}, this equation identifies uniquely
 the pressure field $p$ (as a distribution) up to trivial
 modifications, i.e. additive perturbations depending on time only.

 In the Eulerian-Lagrangian model, instead, the pressure field is
 defined (see (2.20) in \cite{brenier4}) and uniquely determined,
 still up to trivial modifications, by
 \begin{equation}\label{befana2}
 \nabla p(t,x)=-\partial_t\left(\int_D\vv(t,x,a)\,dc_{t,x}(a)\right)
 -\div\left(\int_D\vv(t,x,a)\otimes\vv(t,x,a)\,dc_{t,x}(a)\right),
 \end{equation}
 all derivatives being understood in the sense of distributions in
 $(0,1)\times D$ (here $(c,\vv)$ is any optimal pair for the
 Eulerian-Lagrangian model). We used the same letter $p$ to denote
 the pressure field in the two models: indeed, we have seen in the
 proof of Theorem~\ref{thmequivalence} that, writing
 $\eeta=\eeta_a\otimes\mu_D$, the correspondence
 $$\eeta\mapsto (c^{\seeta}_{t,a},\vv^{\seeta}_{t,a})
 \qquad\text{with}\qquad c^{\seeta}_{t,a}:=(e_t)_\#\eeta_a,\,\,
 \vv^{\seeta}_{t,a}c^{\seeta}_{t,a}:=(e_t)_\#(\dot\omega(t)\eeta_a)
 $$
 maps optimal solutions for the first problem into optimal
 solutions for the second one. Since under this correspondence
 \eqref{befana2} reduces to \eqref{befana}, the two pressure fields
 coincide.

 The following crucial regularity result for the pressure field has
 been obtained in \cite{amfiga}, improving in the time variable
 the regularity $\partial_{x_i} p\in{\cal M}_{\rm loc}\bigl((0,1)\times D\bigr)$
 obtained by Brenier in \cite{brenier4}.

 \begin{theorem}[Regularity of pressure]\label{tregpres}
 Let $(c,\vv)$ be an optimal pair for the Eulerian-La\-gran\-gi\-an
 model, and let $p$ be the pressure field identified by
 \eqref{befana2}. Then $\partial_{x_i}p\in L^2_{\rm loc}\left((0,1);{\cal M}(D)\right)$
 and
 $$
 p\in L^2_{\rm loc}\bigl((0,1);BV_{\rm loc}(D)\bigr)
 \subset
 L^2_{\rm loc}\bigl((0,1);L^{d/(d-1)}_{\rm loc}(D)\bigr).
 $$
 In the case $D=\T^d$ the same properties hold globally in space, i.e.
 replacing $BV_{\rm loc}(D)$ with $BV(\T^d)$ and $L^{d/(d-1)}_{\rm loc}(D)$
 with $L^{d/(d-1)}(\T^d)$.
 \end{theorem}

 The $L^1_{\rm loc}$ integrability of $p$ allows much stronger
 variations in the Lagrangian model, that give rise to possibly
 nonsmooth densities, which may even vanish.

 From now one we shall confine our discussion to the case of the
 flat torus $\T^d$, as our arguments involve some \emph{global}
 smoothing that becomes more technical, and needs to be carefully
 checked in more general situations. We also set
 $\mu_{\T}=\mu_{\T^d}$ and denote by $d_{\T}$ the Riemannian
 distance in $\T^d$ (i.e. the distance modulo 1 in $\R^d/\Z^d$). In
 the next theorem we consider generalized flows $\nnu$ with
 \emph{bounded compression}, defined by the property $\rho^\snnu\in
 L^\infty\left((0,1)\times D\right)$.

 \begin{theorem}
 \label{generalminimality} Let $\eeta$ be an optimal incompressible
 flow in $\T^d$ between $\eta$ and $\gamma$. Then
 \begin{equation}
 \label{eqoptimality1} \langle p,\rho^\snnu - 1 \rangle \leq
 \AA_1(\nnu) - \AA_1(\eeta)
 \end{equation}
 for any generalized flow with bounded compression $\nnu$ between
 $\eta$ and $\gamma$ such that
 \begin{equation}\label{extra}
 \text{$\rho^\snnu(t,\cdot)=1$ for $t$ sufficiently close to $0$,
 $1$.}
 \end{equation}
 If $p\in L^1([0,1]\times\T^d)$, the condition \eqref{extra} is not
 required for the validity of \eqref{eqoptimality1}.
 \end{theorem}
 \begin{Proof}
 Let $J:=\{\rho^{\snnu}(t,\cdot)\neq 1\}\Subset (0,1)$ and let us
 first assume that $\rho^\snnu$ is smooth. If
 $\Vert\rho^\snnu-1\Vert_{C^1}\leq 1/2$, then the result follows by
 Theorem~\ref{thmslighcomprflows}. If not, for $\e>0$ small enough
 $(1-\e)\eeta + \e\nnu$ is a slightly compressible generalized flow
 in the sense of Definition~\ref{defslighcomprflows}. Thus, we have
 $$
 \e \langle p,\rho^\snnu - 1 \rangle= \langle
 p,\rho^{(1-\e)\seeta+\e\snnu} - 1 \rangle \leq
 \AA_1((1-\e)\eeta+\e\nnu) - \AA_1(\eeta) =\e\left(\AA_1(\nnu) -
 \AA_1(\eeta) \right),
 $$
 and this proves the statement whenever $\rho^\snnu$ is smooth.

 If $\rho^\snnu$ is not smooth, we need a regularization argument.
 Let us assume first that $\rho^\snnu$ is smooth in time, uniformly
 with respect to $x$, but not in space. We fix a cut-off function
 $\chi\in C^1_c(0,1)$ strictly positive on a neighbourhood
 of $J$ and define, for $y\in\R^d$, the maps
 $T_{\e,y}:\Oa(\T^d)\to\Oa(\T^d)$ by
 $$
 T_{\e,y}(\omega,a):=(\omega+\e y\chi,a),\quad\qquad
 (\omega,a)\in\O(\T^d).
 $$
 Then, we set $\nnu_{\e}:=\int_{\R^d}(T_{\e,y})_\#\nnu\phi(y)\,dy$,
 where $\phi:\R^d\to [0,+\infty)$ is a standard convolution kernel.
 It is easy to check that $\nnu_\e$ still connects $\eta$ to $\gamma$,
 and that
 $$
 \rho^{\snnu_\e}(t,\cdot)=\rho^{\snnu}(t,\cdot)\ast\phi_{\e\chi(t)}
 \qquad\forall t\in [0,1],
 $$
 where $\phi_\e(x)=\e^{-d}\phi(x/\e)$. Since
 $$
 \lim_{\e\downarrow 0}\AA_1(\nnu_\e)= \lim_{\e\downarrow
 0}\int_{\R^d}\int_{\Oa(\T^d)}\int_0^1 |\dot \omega(t)+\e
 y\dot\chi(t)|^2\,dt\,d\nnu(\omega,a)\phi(y)\,dy =\AA_1(\nnu)
 $$
 we can pass to the limit in \eqref{eqoptimality1} with $\nnu_\e$
 in place of $\nnu$, which are smooth.

 In the general case we fix a convolution kernel with compact
 support $\varphi(t)$ and, with the same choice of $\chi$ done
 before, we define the maps
 $$
 T_{\e}(\omega,a)(t):=(\int_0^1\omega(t-s\e\chi(t))\varphi(s)\,ds,a).
 $$
 Setting $\nnu_\e=(T_\e)_\#\nnu$, it is easy to check that
 $\AA_1(\nnu_\e)\to\AA_1(\nnu)$ and that
 $$
 \rho^{\snnu_\e}(t,x)=\int_0^1\rho^{\snnu}(t-s\e\chi(t),x)\varphi(s)\,ds
 $$
 are smooth in time, uniformly in $x$. So, by applying
 \eqref{eqoptimality1} with $\nnu_\e$ in place of $\nnu$, we obtain
 the inequality in the limit.

 Finally, if $p$ is globally integrable, we can approximate any
 generalized flow with bounded compression $\nnu$ between $\eta$
 and $\gamma$ by transforming $\omega$ into $\omega\circ\psi_\e$,
 where $\psi_\e:[0,1]\to [0,1]$ is defined by
 $\psi_\e(t):=\frac{1}{1-2\e}\int_0^t \chi_{[\e,1-\e]}(s)ds$
 (so that $\psi_\e$ is constant for $t$ close to $0$ and $1$). Passing to
 the limit as $\e\downarrow 0$ we obtain the inequality even
 without the condition $\rho^{\snnu}(t,\cdot)=1$ for $t$ close to
 $0$, $1$.
 \end{Proof}

 \begin{remark}[Smoothing of flows and plans]\label{notsobad}
 {\rm Notice that the same smoothing argument can be used to prove
 this statement: given a flow $\eeta$ between $\eta=\eta_a\otimes\mu_\T$ and
 $\gamma=\gamma_a\otimes\mu_\T$ (not necessarily with bounded
 compression), we can find flows with bounded
 compression $\eeta^\e$ connecting
 $\eta^\e:=(\eta_a)\ast\phi_\e\otimes\mu_\T$ to
 $\gamma^\e:=(\gamma_a)\ast\phi_\e\otimes\mu_\T$, with
 $\AA_T(\eeta^\e)=\AA_T(\eeta)$ and
 $$
 \int_{\Oa(\T^d)}\int_0^1r_\e(\tau,\omega)\,d\tau\,d\eeta(\omega,a)=
 \int_{\Oa(\T^d)}\int_0^1r(\tau,\omega)\,d\tau\,d\eeta^\e(\omega,a)
 \qquad \forall r\in L^1\left([0,1]\times\T^d\right)
 $$
 (where, as usual, $r_\e(t,x)=r(t,\cdot)\ast\phi_\e(x)$). In order
 to have these properties, it suffices to define
 $$
 \eeta^\e:=\int_{\R^d}(\sigma_{\e y})_\#\eeta\,\phi(y)\,dy,
 $$
 where $\sigma_z(\omega,a)=(\omega+z,a)$. Notice also that the
 ``mollified plans'' $\eta^\e,\,\gamma^\e$ converge to $\eta$,
 $\gamma$ in $(\Gamma(\T^d),\overline\delta)$: if we consider the map
 $S^\e_y:\T^d\to\O(\T^d)$ given by $x\mapsto\omega_x(t):=x+\e t y$,
 the generalized incompressible flow
 $\nnu^\e=\nnu^\e_a\otimes\mu_\T$, with
 $$
 \nnu^\e_a:=\int_{\R^d} (S^\e_y)_\#\lambda_a\,\phi(y)\,dy,
 $$
 connects in $[0,1]$ the plan $\lambda=\lambda_a\otimes\mu_\T$ to
 $\lambda^\e=(\lambda_a\ast\phi_\e)\otimes\mu_\T$, with an action
 equal to $\e^2\int_{\R^d}|y|^2\phi(y)\,dy$. }\end{remark}

 In order to state necessary and sufficient optimality conditions
 at the level of single fluid paths, we have to take into account
 that the pressure field is not pointwise defined, and to choose a
 particular representative in its equivalence class, modulo
 negligible sets in spacetime. Henceforth, we define
 \begin{equation}\label{defbarp}
 \bar p(t,x):=\liminf_{\e\downarrow 0}p_\e(t,x),
 \end{equation}
 where, thinking of $p(t,\cdot)$ as a $1$-periodic function in
 $\R^d$, $p_\e$ is defined by
 $$
 p_\e(t,x):=(2\pi)^{-d/2}\int_{\R^d}p(t,x+\e y)e^{-|y|^2/2}\,dy.
 $$
Notice that $p_\e$ is smooth and still 1-periodic. The choice of
the heat kernel here is convenient, because of the semigroup
property $p_{\e+\e'}=(p_\e)_{\e'}$. Recall that $\bar p$ is a
 representative, because at any Lebesgue point $x$ of $p(t,\cdot)$
 the limit of $p_\e(t,x)$ exists, and coincides with $p(t,x)$.

 In order to handle passages to limits, we need also uniform
 pointwise bounds on $p_\e$; therefore we define
 \begin{equation}\label{maximal}
 M f(x):=\sup_{\e>0}\,\,(2\pi)^{-d/2}\int_{\R^d}|f|(x+\e y) e^{-|y|^2/2}\,dy,\quad
 \qquad f\in L^1(\T^d).
 \end{equation}
 We will use the following facts: first, $$
 M f_\e=\sup_{\e'>0}|f_\e|_{\e'}\leq
 \sup_{\e'>0}(|f|_\e)_{\e'}\leq \sup_{r>0}|f|_r=Mf
$$
 because of the semigroup property; second, standard maximal inequalities
 imply $\Vert Mf\Vert_{L^p(\T^d)}\leq c_p\Vert f\Vert_{L^p(\T^d)}$
 for all $p>1$. Setting $Mp(t,x):=Mp(t,\cdot)(x)$, by Theorem \ref{tregpres} we infer that
 $Mp \in L^2_{\rm loc}\bigl((0,1),L^{d/d-1}(\T^d)\bigr)$, so that in particular
 $Mp \in L^1_{\rm loc}\bigl((0,1) \times \T^d\bigr)$. This is the
 integrability assumption on $p$ that will play a role in the rest
 of this section.

 \begin{definition} [$q$-minimizing path]
 Let $\omega\in H^1\left((0,1);D\right)$ with
 $Mq(\tau,\omega)\in L^1(0,1)$.\\
 We say that $\omega$ is a \emph{$q$-minimizing path} if
 $$
 \int_0^1 \frac{1}{2}|\dot\omega(\tau)|^2-q(\tau,\omega)\,d\tau\leq \int_0^1
 \frac{1}{2}|\dot\omega(\tau)+\dot\delta(\tau)|^2-
 q(\tau,\omega+\delta)\,d\tau
 $$
 for all $\delta\in H^1_0\left((0,1);D\right)$
 with $Mq(\tau,\omega+\delta)\in L^1(0,1)$.\\
 Analogously, we say that $\omega$ is a \emph{locally
 $q$-minimizing path} if
 \begin{equation}\label{caffe}
 \int_s^t \frac{1}{2}|\dot\omega(\tau)|^2-q(\tau,\omega)\,d\tau\leq \int_s^t
 \frac{1}{2}|\dot\omega(\tau)+\dot\delta(\tau)|^2-
 q(\tau,\omega+\delta)\,d\tau
 \end{equation}
 for all $[s,t]\subset (0,1)$ and all $\delta\in
 H^1_0\left((s,t);D\right)$ with $Mq(\tau,\omega+\delta)\in
 L^1(s,t)$.
 \end{definition}

 \begin{remark}\label{generic}
 {\rm We notice that, for incompressible flows $\eeta$, the $L^1$ (resp. $L^1_{\rm loc}$)
 integrability of $Mq(\tau,\omega)$ imposed on the curves
 $\omega$ (and on their perturbations $\omega+\delta$) is satisfied $\eeta$-a.e. if
 $Mq\in L^1\left((0,1)\times \T^d\right)$)
 (resp. $Mq\in L^1_{\rm loc}\left((0,1)\times \T^d\right)$);
 this can simply be obtained first noticing that
 the incompressibility of $\eeta$ and Fubini's theorem give
 $$
 \int_{\Oa(\T^d)}\int_Jf(\tau,\omega)\,d\tau\,d\eeta(\omega,a)=
 \int_J\int_{\T^d} f(\tau,x)\,d\mu_{\T^d}(x)\,d\tau
 $$
 for all nonnegative Borel functions $f$ and all intervals
 $J\subset (0,1)$, and then applying this identity to $f=Mq$.}
 \end{remark}

 \begin{theorem}[First necessary condition]\label{optima1}
 Let $\eeta=\eeta_a\otimes\mu_\T$ be any optimal incompressible
 flow on $\T^d$. Then, $\eeta$ is concentrated on locally $\bar
 p$-minimizing paths, where $\bar p$ is the precise representative
 of the pressure field $p$, and on $\bar p$-minimizing paths if
 $Mp\in L^1([0,1]\times\T^d)$.
 \end{theorem}
 \begin{Proof} With no loss of generality we identify $\T^d$ with $\R^d/\Z^d$.
 Let $\eeta$ be an optimal incompressible flow and $[s,t]\subset
 (0,1)$. We fix a nonnegative function $\chi\in C^1_c(0,1)$ with
 $\{\chi >0\}=(s,t)$. Given
 $\delta\in H^1_0\left([s,t];\T^d\right)$, $y\in\R^d$ and a Borel set
 $E\subset\Oa(\T^d)$, we define $T_{\e,y}:\Oa(\T^d)\to\Oa(\T^d)$ by
 $$
 T_{\e,y}(\omega,a):=
 \begin{cases}
 (\omega,a)&\text{if $\omega\notin E$;}\\
 (\omega+\delta+\e y\chi,a) &\text{if $\omega\in E$}
 \end{cases}
 $$
 (of course, the sum is understood modulo $1$) and
 $\nnu_{\e,y}:=(T_{\e,y})_\#\eeta$.

 It is easy to see that $\nnu_{\e,y}$ is a flow with bounded
 compression, since for all times $\tau$ the curves $\omega(\tau)$
 are either left unchanged, or translated by the constant
 $\delta(\tau)+\e y\chi(\tau)$, so that the density produced by
 $\nnu_{\e,y}$ is at most 2, and equal to $1$ outside the interval
 $[s,t]$.

 Therefore, by Theorem~\ref{generalminimality} we get
 $$
 \int_{\T^d}\int_s^t\bar p
 (\rho^{\snnu_{\e,y}}-1)\,d\tau\,d\mu_\T\leq \int_E \AA_1(\omega +
 \delta+\e y\chi)-\AA_1(\omega) \,d\eeta(\omega,a).
 $$
 Rearranging terms, we get
 $$
 \int_E \int_s^t\frac{1}{2}|\dot\omega|^2-\bar
 p(\tau,\omega)\,d\tau\,d\eeta(\omega,a) \leq\int_E
 \left[\int_s^t\frac{1}{2}|\dot\omega+\dot\delta+\e y\dot\chi|^2- \bar
 p(\tau,\omega+\delta+\e y\chi)\,d\tau\right]d\eeta(\omega,a).
 $$
 We can now average the above inequality using the
 heat kernel $\phi(y)=(2\pi)^{-d/2}e^{-|y|^2/2}$, and
 we obtain
 \begin{multline*}
 \int_E \int_s^t\frac{1}{2}|\dot\omega|^2-\bar p(\tau,\omega)\,d\tau\,d\eeta(\omega,a)\\
 \leq \int_E \left[\int_{\R^d}\int_s^t\frac{1}{2}|\dot\omega+\dot\delta+\e
 y\dot\chi|^2\,d\tau \phi(y)\,dy -\int_s^t
 p_{\e\chi(\tau)}(\tau,\omega+\delta)\,d\tau\right]
 d\eeta(\omega,a).
 \end{multline*}
 Now, let $\mathcal D\subset H^1_0\left([s,t];\T^d\right)$ be a
 countable dense subset; by the arbitrariness of $E$ and
 Remark~\ref{generic} we infer the existence of a
 $\eeta$-negligible Borel set $B\subset\Oa(\T^d)$ such that
 $Mp(\tau,\omega)\in L^1(s,t)$ and
 $$
 \int_s^t\frac{1}{2}|\dot\omega|^2-\bar p(\tau,\omega)\,d\tau\leq
 \int_{\R^d}\int_s^t \frac{1}{2}|\dot\omega+\dot\delta+\e y\dot\chi|^2\,d\tau
 \phi(y)\,dy-\int_s^t p_{\e\chi(\tau)}(\tau,\omega+\delta)\,d\tau
 $$
 holds for all $\e=1/n$, $\delta\in\mathcal D$ and
 $(\omega,a)\in\Oa(\T^d)\setminus B$. By a density argument, we see
 that the same inequality holds for all $\e=1/n$, $\delta\in
 H^1_0\left([s,t];\T^d\right)$, and $\omega\in\Oa(\T^d)\setminus
 B$.

 Now, if $Mp(\tau,\omega+\delta)\in L^1(s,t)$, since $\delta\in
 H^1_0\left([s,t];\T^d\right)$ we have that
 $Mp(\tau,\omega+\delta)\in L^1(s,t)$, and we can use the bound
 $|p_\e|\leq Mp$ to pass to the limit as $\e\downarrow 0$
  to obtain that \eqref{caffe} holds
 with $q=\bar p$.

 The proof of the global minimality property in the case when $p\in
 L^1([0,T]\times\T^d)$ is similar, just letting $\delta$ vary in
 $H^1_0\left([0,1];\T^d\right)$ and using a fixed function $\chi\in
 C^1([0,1])$ with $\chi(0)=\chi(1)=0$ and $\chi>0$ in $(0,1)$.
 \end{Proof}

 In order to state the second necessary optimality condition
 fulfilled by minimizers, we need some preliminary definition.
 Let $q\in L^1\left([s,t]\times D\right)$ and let us define the cost
 $c^{s,t}_q:D\times D\to\overline\R$ of the minimal connection in
 $[s,t]$ between $x$ and $y$, namely
 \begin{equation}\label{olla}
 c^{s,t}_q(x,y):=\inf\left\{
 \int_s^t\frac{1}{2}|\dot\omega(\tau)|^2-q(\tau,\omega)\,d\tau:\
 \omega(s)=x,\,\omega(t)=y,\,Mq(\tau,\omega)\in L^1(s,t)\right\},
 \end{equation}
 with the convention $c^{s,t}_q(x,y)=+\infty$ if no admissible
 curve $\omega$ exists. Using this cost function
 $c^{s,t}_q$, we can consider the induced
 optimal transport problem, namely
 \begin{equation}\label{defwc}
 W_{c^{s,t}_q}(\mu_1,\mu_2):=\inf\left\{\int_{D\times D}
 c_q^{s,t}(x,y)\,d\lambda(x,y):\ \lambda\in\Gamma(\mu_1,\mu_2),
 \,\,(c_q^{s,t})^+\in L^1(\lambda)\right\},
 \end{equation}
 where $\Gamma(\mu_1,\mu_2)$ is the family of all probability
 measures $\lambda$ in $D\times D$ whose first and second marginals
 are respectively $\mu_1$ and $\mu_2$. Again, we set by convention
 $W_{c^{s,t}_q}(\mu_1,\mu_2)=+\infty$ if no admissible $\lambda$
 exists.

 Unlike most classical situations (see \cite{villani}), existence
 of an optimal $\lambda$ is not guaranteed because $c^{s,t}_q$ are
 not lower semicontinuous in $D\times D$, and also it seems difficult
 to get lower bounds on $c_q^{s,t}$. It will be useful,
 however, the following \emph{upper} bound on $W_{c^{s,t}_q}$:

 \begin{lemma} \label{lspread}
 If $Mq\in L^1([s,t]\times\T^d)$ there exists a nonnegative $\mu_\T$-integrable
 function $K^{s,t}_q$ satisfying
 \begin{equation}\label{dominationcost}
 c^{s,t}_q(x,y) \leq K_q^{s,t}(x) + K_q^{s,t}(y) \qquad\forall x,y \in \T^d.
 \end{equation}
 \end{lemma}

 \begin{remark}\label{rmkdomination}
 {\rm By \eqref{dominationcost} we deduce that, if $K^{s,t}_q\in L^1(\mu_1+\mu_2)$, then
 $(c_q^{s,t})^+\in L^1(\lambda)$ for all $\lambda\in\Gamma(\mu_1,\mu_2)$ and we have
 $$
 \int_{\T^d\times\T^d}c^{s,t}_q(x,y)\,d\lambda(x,y)\leq
 \int_{\T^d}K_q^{s,t}(w)\,d(\mu_1+\mu_2)(w) \qquad\forall
 \lambda\in\Gamma(\mu_1,\mu_2).
 $$
 In particular, $W_{c^{s,t}_q}(\mu_1,\mu_2)$ as defined in \eqref{defwc} is not
 equal to $+\infty$.}
 \end{remark}

 \begin{Proof}
 Assume $s=0$ and let $l=t/2$. Let us fix
 $x,\,y\in\T^d$; given $z\in\T^d$ we consider the projection
 on $\T^d$ of the Euclidean path
 $$
 \omega_z(\tau):=
 \begin{cases}
 x+\frac{\tau}{l}(z-x)&\text{if $\tau\in [0,l]$;}\\
 z+\frac{\tau-l}{l}(y-z)&\text{if $\tau\in [l,t]$.}
 \end{cases}
 $$
 This path leads to the estimate
 $$
 c^{0,t}_q(x,y)\leq\frac{d_\T^2(x,z)+d_\T^2(z,y)}{2l}+
 \int_0^lMq(\tau,x+\frac{\tau}{l}(z-x))\,d\tau+
 \int_l^tMq(\tau,z+\frac{\tau-l}{l}(y-z))\,d\tau.
 $$
 By integrating the free variable $z$ with respect to $\mu_\T$,
 since $d_\T\leq\frac{\sqrt{d}}{2}$ on $\T^d\times\T^d$, we get
 $$
 c^{0,t}_q(x,y)\leq \frac{d}{4l}+
 \int_{\T^d}\int_0^lMq(\tau,x+\frac{\tau}{l}(z-x))+
 Mq(l+\tau,z+\frac{\tau}{l}(y-z))\,d\tau\,d\mu_\T(z).
 $$
 Therefore, the function
 \begin{equation}\label{defkq}
 K^{0,t}_q(w):=\frac{d}{4l}+\int_{\T^d}\int_0^lMq(\tau,w+\frac{\tau}{l}(z-w))+
 Mq(l+\tau,z+\frac{\tau}{l}(w-z))\,d\tau\,d\mu_\T(z)
 \end{equation}
 fulfils \eqref{dominationcost}. It is easy to check, using Fubini's
 theorem, that $K^{0,t}_q$ is $\mu_\T$-integrable in $\T^d$.
 Indeed,
 \begin{align*}
 \int_{\T^d} K^{0,t}_q(w)\,d\mu_\T(w)
 &=\frac{d}{4l}+\int_{\T^d}\int_{\T^d}\int_0^l
 Mq(\tau,w+\frac{\tau}{l}(z-w))\,d\tau\,d\mu_\T(z)\,d\mu_\T(w)\\
 &+\int_{\T^d}\int_{\T^d}\int_0^l
 Mq(l+\tau,z+\frac{\tau}{l}(w-z))\,d\tau\,d\mu_\T(w)\,d\mu_\T(z)\\
 &=\frac{d}{4l}+\int_{\T^d}\int_{\T^d}\int_0^l
 Mq(\tau,w+\frac{\tau}{l}y)\,d\tau\,d\mu_\T(y)\,d\mu_\T(w)\\
 &+\int_{\T^d}\int_{\T^d}\int_0^l
 Mq(l+\tau,z+\frac{\tau}{l}y)\,d\tau\,d\mu_\T(z)\,d\mu_\T(y)\\
 &=\frac{d}{4l}+\int_0^l\int_{\T^d}\int_{\T^d}
 Mq(\tau,w+\frac{\tau}{l}y)+
 Mq(l+\tau,w+\frac{\tau}{l}y)\,d\mu_\T(w)\,d\mu_\T(y)\,d\tau\\
 &=\frac{d}{4l}+\int_0^t\int_{\T^d} Mq(\tau,w)\,d\mu_\T(w)\,d\tau<+\i.
 \end{align*}
 \end{Proof}

 In the proof of the next theorem we are going to use the
 \emph{measurable selection theorem} (see \cite[Theorems
 III.22 and III.23]{castval}): if $(A,\mathcal A,\nu)$ is a measure space, $X$
 is a Polish space and $E\subset A\times X$ is $\mathcal
 A_\nu\otimes \BorelSets{X}$-measurable, where $\mathcal A_\nu$ is
 the $\nu$-completion of $\mathcal A$, then:
 \begin{itemize}
 \item[(i)] the projection $\pi_A(E)$ of $E$ on $A$
 is $\mathcal A_\nu$-measurable;
 \item[(ii)] there exists a $(\mathcal
 A_\nu,\BorelSets{X})$-measurable map $\sigma:\pi(E)\to X$ such
 that $(x,\s(x))\in E$ for $\nu$-a.e. $x\in \pi_A(E)$.
 \end{itemize}

 The next theorem will provide a new necessary optimality
 condition involving not only the path that should be followed
 between $x$ and $y$ (which, as we proved, should minimize the
 Lagrangian ${\cal L}_{\bar p}$ in
 \eqref{defLp}), but also the ``weights'' given to the paths. We observe that,
 when
 a variation of these weights is performed, new flows $\tilde{\eeta}$ between $\eta$ and $\g$
 are built which need not be of bounded compression, for which $(e_t)_\#\tilde{\eeta}$
 might be even singular with respect to $\mu_{\T}$; therefore
 we can't use directly them in the variational principle \eqref{eqoptimality1};
 however, this difficulty can be overcome by the smoothing procedure
 in Remark~\ref{notsobad}.

 \begin{theorem}[Second necessary condition]\label{optima2}
 Let $\eeta=\eeta_a\otimes\mu_\T$ be an optimal incompressible flow
 on $\T^d$ between $\eta$ and $\g$. Then, for all intervals $[s,t]\subset (0,1)$,
 $W_{c^{s,t}_{\bar p}}(\eta_a,\gamma_a)\in\R$ and the plan
 $(e_s,e_t)_\#\eeta_a$ is optimal, relative to the cost
 $c^{s,t}_{\bar p}$ defined in \eqref{olla}, for $\mu_\T$-a.e. $a$.
 \end{theorem}
 \begin{Proof} Let $[s,t]\subset (0,1)$ be fixed. Since
 \begin{eqnarray}\label{chefffatica}
 \int_{\T^d}\int_{\T^d\times\T^d}c_{\bar
 p}^{s,t}(x,y)\,d(e_s,e_t)_\#\eeta_a
 \,d\mu_\T(a)&\leq&\int_{\T^d}\int_{\O(\T^d)}\int_s^t
 \frac{1}{2}|\dot\omega(\tau)|^2-\bar p(\tau,\omega)
 \,d\tau d\eeta_a(\omega)\,d\mu_\T(a)\nonumber \\
 &=&(t-s)\overline\delta^2(\eta,\gamma),
 \end{eqnarray}
 it suffices to show that
 \begin{equation}\label{cheffatica}
 (t-s)\overline\delta^2(\eta,\gamma)\leq\int_{\T^d}
 W_{c^{s,t}_{\bar p}}(\eta^s_a,\gamma^t_a)\,d\mu_\T(a).
 \end{equation}
 We are going to prove this fact by a smoothing argument. We set
 $\eta^s=\eta^s_a\otimes\mu_\T$,
 $\gamma^t=\gamma^t_a\otimes\mu_\T$, with
 $\eta^s_a=(e_s)_\#\eeta_a$, $\gamma^t_a=(e_t)_\#\eeta_a$. Recall
 that Remark~\ref{rdavid} gives
 $$
 \overline\delta(\eta,\eta^s)=s\overline\delta(\eta,\gamma),\qquad
 \overline\delta(\gamma^t,\gamma)=(1-t)\overline\delta(\eta,\gamma).
 $$
 First, we notice that Lemma~\ref{lspread} gives
 \begin{equation}\label{disacco}
 \begin{split}
 \int_{\T^d}W_{c^{s,t}_{-|\bar p|}}(\eta^s_a,\gamma^t_a)\,d\mu_\T(a)&\leq
 \int_{\T^d}\int_{\T^d}K^{s,t}_{-|\bar p|}(w)\,d(\eta^s_a+\gamma^t_a)(w)\,d\mu_\T(a)\\
 &=2\int_{\T^d}K^{s,t}_{-|\bar p|}(w)\,d\mu_\T(w)<+\infty.
 \end{split}
 \end{equation}
 We also remark that, since $\tau\mapsto\Vert p_\e(\tau,\cdot)\Vert_\infty$ is
 integrable in $(s,t)$, for any $\e>0$ the cost $c^{s,t}_{p_\e}$ is bounded both from above
and below.
 Next, we show that
 \begin{equation}\label{disacco1}
 c^{s,t}_{\bar p}(x,y)\geq\limsup_{\e\downarrow 0}
 c^{s,t}_{p_\e}(x,y) \qquad\forall (x,y)\in\T^d\times\T^d.
 \end{equation}
 Indeed, let $\omega\in H^1\left([s,t];\T^d\right)$ with
 $\omega(s)=x$, $\omega(t)=y$ and $Mp(\tau,\omega)\in L^1(s,t)$ (if
 there is no such $\omega$, there is nothing to prove). By the
 pointwise bound $|p_\e|\leq Mp$ and Lebesgue's theorem, we get
 $$
 \int_s^t\frac{1}{2}|\dot\omega(\tau)|^2-\bar p(\tau,\omega)\,d\tau
 =\lim_{\e\downarrow 0}
 \int_s^t\frac{1}{2}|\dot\omega(\tau)|^2- p_\e(\tau,\omega)\,d\tau.
 $$
 By the $L^1(L^\infty)$ bound on $Mp_\e$,
 the curve $\omega$ is admissible also for the
 variational problem defining $c_{p_\e}^{s,t}$, therefore the above
 limit provides an upper bound on $\limsup_\e
 c^{s,t}_{p_\e}(x,y)$. By minimizing with respect to $\omega$ we
 obtain \eqref{disacco1}.

 By \eqref{disacco} and the pointwise bound $\bar p\geq -|\bar p|$ we
 infer that the positive part of $W_{c^{s,t}_{\bar
 p}}(\eta^s_a,\gamma^t_a)$ is $\mu_\T$-integrable. Let now
 $\delta>0$ be fixed, and let us consider the compact space
 $X:=\Probabilities{\T^d\times\T^d}$ and the
 $\BorelSets{\T^d}_{\mu_\T}\otimes\BorelSets{X}$-measurable set
 $$
 E:=\left\{(a,\lambda)\in\T^d\times X:\
 \lambda\in\Gamma(\eta^s_a,\gamma^t_a),
 \,\,\int_{\T^d\times\T^d}c^{s,t}_{\bar p}(x,y)\,d\lambda<
 \delta+ \Bigl( W_{c^{s,t}_{\bar p}}(\eta^s_a,\gamma^t_a)\lor -\frac{1}{\delta} \Bigr) \right\}
 $$
 (we skip the proof of the measurability, that is based on tedious but routine
 arguments). Since $W_{c^{s,t}_{\bar p}}(\eta^s_a,\gamma^t_a)<+\infty$
 for $\mu_\T$-a.e. $a$, we obtain that for $\mu_\T$-a.e. $a\in\T^d$
 there exists $\lambda\in\Gamma(\eta^s_a,\gamma^t_a)$ with
 $(a,\lambda)\in E$. Thanks to the measurable selection theorem we can
 select a Borel family
 $a\mapsto\lambda_a\in\Probabilities{\T^d\times\T^d}$ such that
 $\lambda_a\in\Gamma(\eta^s_a,\gamma^t_a)$ and
 $$
 \int_{\T^d\times\T^d}c^{s,t}_{\bar p}(x,y)\,d\lambda_a<
 \delta+\Bigl( W_{c^{s,t}_{\bar p}}(\eta^s_a,\gamma^t_a)\lor -\frac{1}{\delta} \Bigr)
 \qquad\text{for $\mu_\T$-a.e. $a\in \T^d$.}
 $$
 By Lemma~\ref{lspread} and Remark~\ref{rmkdomination} we get
 $$
 c^{s,t}_{p_\e}(x,y)\leq K_{p_\e}^{s,t}(x) + K_{p_\e}^{s,t}(y)
 \leq K_{p}^{s,t}(x) + K_{p}^{s,t}(y) \qquad \forall x,y \in \T^d
 $$
 and
 $$
 \int_{\T^d}\int_{\T^d\times\T^d}K_{p}^{s,t}(x) + K_{p}^{s,t}(y)\,d\lambda_a\,d\mu_{\T}(a)=
 \int_{\T^d}\int_{\T^d} K_p^{s,t}\,d(\eta^s_a+\gamma^t_a)\,d\mu_\T(a)<+\infty
 $$
 (we used the pointwise bound $Mp_\e\leq Mp$ and the fact that
 $q\mapsto K_q^{s,t}$ has a monotone dependence upon $Mq$, see
 \eqref{defkq}). Therefore
 \eqref{disacco1} and Fatou's lemma give
 \begin{equation}\label{disacco2}
 \delta+\int_{\T^d}\Bigl( W_{c^{s,t}_{\bar p}}(\eta^s_a,\gamma^t_a)\lor -\frac{1}{\delta} \Bigr)
 \,d\mu_\T(a)\geq\limsup_{\e\downarrow 0} \int_{\T^d}
 \int_{\T^d\times\T^d}c_{p_\e}^{s,t}(x,y) d\lambda_a\,d\mu_\T(a).
 \end{equation}
 Still thanks to the measurable selection theorem, we can find a
 Borel map $(x,y,a)\mapsto\omega^{x,y}_{a,\e}\in
 C\left([s,t];\T^d\right)$ with $\omega^{x,y}_{a,\e}(s)=x$,
 $\omega^{x,y}_{a,\e}(t)=y$, $Mp_\e(\tau,\omega^{x,y}_{a,\e})\in
 L^1(s,t)$ and
 $$
 \int_s^t\frac{1}{2}|\dot\omega^{x,y}_{a,\e}|^2
 -p_\e(\tau,\omega^{x,y}_{a,\e})\,d\tau<\delta+c^{s,t}_{p_\e}(x,y)
 \qquad\text{for $\lambda_a\otimes\mu_\T$-a.e. $(x,y,a)$.}
 $$
 Let $\llambda^\e=\llambda^\e_a\otimes\mu_\T$ be the push-forward,
 under the map $(x,y,a)\mapsto\omega^{x,y}_{a,\e}$, of the measure
 $\lambda_a\otimes\mu_\T$; by construction this measure fulfils
 $(e_s)_\#\llambda^\e_a=\eta^s_a$,
 $(e_t)_\#\llambda^\e_a=\gamma^t_a$, (because the marginals of
 $\lambda_a$ are $\eta^s_a$ and $\gamma^t_a$), therefore it
 connects $\eta^s$ to $\gamma^t$ in $[s,t]$. Then, from
 \eqref{disacco2} we get
 $$
 2\delta+\int_{\T^d}\Bigl( W_{c^{s,t}_{\bar p}}(\eta^s_a,\gamma^t_a)
 \lor -\frac{1}{\delta} \Bigr) \,d\mu_\T(a)\geq \limsup_{\e\downarrow 0}
 \int_{C([s,t];\T^d)\times\T^d}
 \int_s^t\frac{1}{2}|\dot\omega(\tau)|^2-p_\e(\tau,\omega)
 \,d\tau\,d\llambda^\e(\omega,a).
 $$
 Eventually, Remark~\ref{notsobad} provides us with a flow with
 bounded compression $\hat{\llambda}^\e$ connecting $\eta^{s,\e}$
 to $\gamma^{t,\e}$ in $[s,t]$ with
 \begin{equation}\label{disacco3}
 2\delta+\int_{\T^d}
 \Bigl( W_{c^{s,t}_{\bar p}}(\eta_a,\gamma_a)\lor -\frac{1}{\delta} \Bigr)\,d\mu_\T(a)
 \geq \limsup_{\e\downarrow 0}\int_{C([s,t];\T^d)\times\T^d}
 \int_s^t\frac{1}{2}|\dot\omega(\tau)|^2-\bar p(\tau,\omega)
 \,d\tau\,d\hat{\llambda}^\e(\omega,a).
 \end{equation}
 Since $\eta^{s,\e}\to\eta^s$ and $\gamma^{t,\e}\to\gamma^t$ in
 $(\Gamma(\T^d),\overline\delta)$, we can find (by scaling $\eeta$
 from $[0,s]$ to $[0,s_\e]$ and from $[t,1]$ to $[t_\e,1]$, and
 using repeatedly the concatenation, see Remark \ref{rlocalgeo})
 generalized flows $\nnu^\e$ between $\gamma$
 and $\eta$ in $[0,1]$, $s_\e\uparrow s$, $t_\e\downarrow t$
 satisfying:
 \begin{itemize}
 \item[(a)] $\nnu^\e$ connects $\eta$ to $\eta^s$ in $[0,s_\e]$,
 $\eta^s$ to $\eta^{s,\e}$ in $[s_\e,s]$, $\gamma^{t,\e}$ to
 $\gamma^t$ in $[t,t_\e]$, $\gamma^t$ to $\gamma$ in $[t_\e,1]$ and
 is incompressible in all these time intervals;
 \item[(b)] the
 restriction of $\nnu^\e$ to $[s,t]$ coincides with
 $\hat{\llambda}^\e$;
 \item[(c)] the action of $\nnu^\e$ in $[0,s]$
 converges to
 $\overline\delta^2(\eta,\eta^s)=s^2\overline\delta^2(\eta,\gamma)$,
 and the action of $\nnu^\e$ in $[t,1]$ converges to
 $\overline\delta^2(\gamma^t,\gamma)=(1-t)^2\overline\delta^2(\eta,\gamma)$.
 \end{itemize}
 Since $\nnu^\e$ is a flow with bounded compression connecting
 $\eta$ to $\gamma$ we use \eqref{eqoptimality1} and the
 incompressibility in $[0,1]\setminus [s,t]$ to obtain
 \begin{equation}\label{disacco4}
 \int_{\Oa(\T^d)}\int_0^1\frac{1}{2}|\dot\omega(\tau)|^2\,d\tau\,d\nnu^\e(\omega,a)-
 \int_{\Oa(\T^d)}\int_s^t \bar p(\tau,\omega)\,d\tau\,d\nnu^\e(\omega,a)\geq
 \overline\delta^2(\eta,\gamma)
 \end{equation}
 for all $\e>0$.
 Taking into account that (b) and (c) imply
 $$
 \int_{\Oa(\T^d)}\int_0^s \frac{1}{2} |\dot\o(\tau)|^2\,d\tau\,d\nnu^\e(\omega,a) \to s\overline\delta^2(\eta,\gamma)
 $$
 and
 $$
 \int_{\Oa(\T^d)}\int_t^1 \frac{1}{2} |\dot\o(\tau)|^2\,d\tau\,d\nnu^\e(\omega,a) \to (1-t)\overline\delta^2(\eta,\gamma),
 $$
 from \eqref{disacco3} and \eqref{disacco4} we get
 \begin{equation}\label{disacco5}
 2\delta+\int_{\T^d}\Bigl( W_{c^{s,t}_{\bar p}}(\eta^s_a,\gamma^t_a)\lor -\frac{1}{\delta} \Bigr)
 \,d\mu_\T(a) \geq (1-s-(1-t))
 \overline\delta^2(\eta,\gamma)=
 (t-s)\overline\delta^2(\eta,\gamma).
 \end{equation}
 Letting $\delta\downarrow 0$ we obtain the $\mu_\T$-integrability of
 $W_{c^{s,t}_{\bar p}}(\eta^s_a,\gamma^t_a)$ and
 \eqref{cheffatica}.
 \end{Proof}

 A byproduct of the above proof is that equalities hold in
 \eqref{chefffatica}, \eqref{cheffatica}, and therefore
 \begin{eqnarray}\label{insu}
 &&\int_{\T^d}\int_{\O(\T^d)}\left(\int_s^t\frac{1}{2}|\dot\omega(\tau)|^2
 -\bar p(\tau,\omega)\,d\tau-c^{s,t}_{\bar p}(\omega(s),\omega(t))\right)
 d\eeta_a(\omega)\,d\mu_\T(a)\\&&
 =\int_{\T^d}\int_{\O(\T^d)}\int_s^t\frac{1}{2}|\dot\omega(\tau)|^2
 -\bar p(\tau,\omega)\,d\tau\,d\eeta_a(\omega)\,d\mu_\T(a)-
 \int_{\T^d} W_{c^{s,t}_{\bar p}}(\eta^s_a,\gamma^t_a)\,d\mu_\T(a)=0.\nonumber
 \end{eqnarray}
 This yields in particular also the first optimality condition.
 However, as the proof of Theorem~\ref{optima2} is much more
 technical than the one presented in Theorem~\ref{optima1}, we
 decided to present both.

 Now we show that the optimality conditions in
 Theorems~\ref{optima1} and \ref{optima2} are also
 sufficient, even in the case of a general compact manifold without
 boundary $D$.

 \begin{theorem}[Sufficient condition]\label{optima3}
 Assume that $\eeta=\eeta_a\otimes\mu$ is a generalized
 incompressible flow in $D$ between $\eta$ and $\gamma$, and assume
 that for some map $q$ the following properties hold:
 \begin{itemize}
 \item[(a)] $Mq\in L^1((0,1)\times D)$ and $\eeta$ is concentrated
 on $q$-minimizing paths;
 \item[(b)] the plan $(e_0,e_1)_\#\eeta_a$
 is optimal, relative to the cost $c^{0,1}_q$ defined in
 \eqref{olla}, for $\mu_D$-a.e. $a$.
 \end{itemize}
 Then $\eeta$ is optimal and $q$ is the pressure field. In addition,
 if (a), (b) are replaced by
 \begin{itemize}
 \item[(a')] $Mq\in L^1_{\rm loc}((0,1)\times D)$ and $\eeta$ is
 concentrated on locally $q$-minimizing paths;
 \item[(b')] for all
 intervals $[s,t]\subset (0,1)$, the plan $(e_s,e_t)_\#\eeta_a$ is
 optimal, relative to the cost $c^{s,t}_q$ defined in \eqref{olla},
 for $\mu_D$-a.e. $a$,
 \end{itemize}
 the same conclusions hold.
 \end{theorem}
 \begin{Proof}
 Assume first that (a) and (b) hold, and assume without loss of
 generality that $\int_D q(t,\cdot)\,d\mu_D=0$ for almost all $t\in
 (0,1)$. Recalling that, thanks to the global integrability of $Mq$,
 any generalized incompressible flow $\nnu=\nnu_a\otimes\mu_D$
 between $\eta$ and $\gamma$ is concentrated on curves $\omega$
 with $Mq(\tau,\omega)\in L^1(0,1)$ (see Remark~\ref{generic}), we have
 \begin{eqnarray}\label{insu1}
 \AA_1(\nnu)&=&
 \int_D\int_{\O(D)}\int_0^1\frac{1}{2}|\dot\omega|^2-
 q(\tau,\omega)\,d\tau\,d\nnu_a(\omega)\,d\mu_D(a)
 \\
 &\geq& \int_D \int_{D\times
 D}c_q^{0,1}(x,y)\,d(e_0,e_1)_\#\nnu_a\,d\mu_D(a) \geq\int_D
 W_{c_q^{0,1}}(\eta_a,\gamma_a)\,d\mu_D(a).\nonumber
 \end{eqnarray}
 When $\nnu=\eeta$ the first inequality is an equality, because
 $\eeta$ is concentrated on $q$-minimizing paths, as well as the
 second inequality, because of the optimality of the plan
 $(e_0,e_1)_\#\eeta_a$. This proves that $\eeta$ is optimal.
 Moreover, by using the inequality in \eqref{insu1} with a flow $\nnu$
 with bounded compression, one obtains
 $$
 \AA_1(\nnu)\geq\AA_1(\eeta)+\langle q,\rho^{\snnu}-1\rangle.
 $$
 Considering almost incompressible flows $\nnu$ arising by a smooth
 perturbation of $\eeta$ as described at the beginning of this
 section (see \eqref{variazw} in particular), the same
 argument used to obtain \eqref{befana} gives that $q$ satisfies
 \eqref{befana}, so that $q$ is the pressure field.

 In the case when (a)' and (b)' hold, by localizing in all
 intervals $[s,t]\subset (0,1)$ the previous argument (see
 Remark~\ref{rlocalgeo}), one obtains that
 $$
 (t-s)\int_{\Oa(D)}\int_s^t\frac{1}{2}|\dot\omega|^2\,d\tau\,d\eeta(\omega,a)=
 \overline\delta^2(\gamma_s,\gamma_t),
 $$
 where $\gamma_s=(e_s,\pi_D)_\#\eeta$ and
 $\gamma_t=(e_t,\pi_D)_\#\eeta$. Letting $s\downarrow 0$ and
 $t\uparrow 1$ we obtain the optimality of $\eeta$.
 \end{Proof}

 A byproduct of the previous result is a new variational principle
 satisfied, at least locally in time, by the pressure field. Up to a restriction
 to a smaller time interval we shall assume that $Mp\in L^1\left([0,1]\times\T^d\right)$.

 \begin{corollary}[Variational characterization of the pressure]\label{vpre}
 Let $\eta,\,\gamma\in\Gamma(\T^d)$ and let $p$ be the unique pressure field induced by
 the constant speed geodesics in $[0,1]$ between $\eta=\eta_a\otimes\mu_\T$ and
 $\gamma=\gamma_a\otimes\mu_\T$. Assume
 that $Mp\in L^1\left([0,1]\times\T^d\right)$ and, with no loss of generality,
 $\int_{\T^d}p(t,\cdot)\,d\mu_\T=0$. Then
 $\bar p$ maximizes the functional
 $$
 q\mapsto \Psi(q):=\int_{\T^d}W_{c^{0,1}_q}(\eta_a,\gamma_a)\,d\mu_\T(a)+
 \int_0^1\int_{\T^d}q(\tau,x)\,d\mu_\T(x)\,d\tau
 $$
 among all functions $q:[0,1]\times\T^d\to\R$ with $Mq\in L^1([0,1]\times\T^d)$.
 \end{corollary}
 \begin{Proof}
 We first remark that the functional $\Psi$ is
 invariant under sum of functions depending on $t$ only,
 so we can assume that the spatial means of any function $q$ vanish.\\
 From \eqref{insu} we obtain that
 $$
 \int_{\T^d}W_{c^{0,1}_{\bar p}}(\eta_a,\gamma_a)\,d\mu_\T(a)
 =\int_{\T^d}\int_{\O(\T^d)}\int_0^1\frac{1}{2}|\dot\omega(\tau)|^2
 -\bar p(\tau,\omega)\,d\tau\,d\eeta_a(\omega)\,d\mu_\T(a).
 $$
 By the incompressibility constraint, in the right hand side $\bar p$ can be replaced by
 any function $q$ whose spatial means vanish and, if $Mq\in L^1\left([0,1]\times\T^d\right)$,
 the resulting integral bounds from above
 $\int_{\T^d}W_{c^{0,1}_q}(\eta_a,\gamma_a)\,d\mu_\T(a)$,
 as we proved in \eqref{insu1}.
 \end{Proof}

 \end{document}